\documentclass[11pt]{amsart}
\usepackage{amsthm}
\usepackage{amssymb}
\usepackage{amsmath}
\usepackage{graphicx}
\usepackage{amscd}
\usepackage{url}

\newcommand     {\comment}[1]   {}
\newcommand{\mute}[2] {}
\newcommand     {\printname}[1] {}
\newcommand{\labell}[1] {\label{#1}\printname{#1}}

\numberwithin{equation}{section}
\newtheorem {Theorem}[equation]         {Theorem}

\newtheorem {Lemma}[equation]           {Lemma}

\newtheorem {Corollary} [equation]      {Corollary}

\newtheorem* {Theorem*}                  {Theorem}

\theoremstyle{definition}

\theoremstyle{remark}
\newtheorem{Remark}[equation]{Remark}

\newtheorem{Example}[equation]{Example}

\def \del {\partial}
\newcommand {\dd}[1] {\frac{\del}{\del #1}}

\def    \real {{\operatorname{Re}}}

\def	\Flux   {\operatorname{Flux}}
\def	\Ham	{\operatorname{Ham}}

\def	\length	{\operatorname{length}}

\def    \rank   {\operatorname{rank}}
\def    \codim  {\operatorname{codim}}
\def    \calN   {{\mathcal N}}
\def	\Bmax	{B_{\max}}
\def	\Nmax	{\Bmax}
\def	\Neps	{{\mathcal N}_{\eps}}
\def    \pmax   {p_{\text{max}}}

\def    \epsbar {\ol{\eps}}
\def    \phibar {\ol{\phi}}
\def    \psibar {\ol{\psi}}

\def    \g      {{\mathfrak g}}
\def    \hatN   {{\hat{N}}}

\def	\Z	{\mathbb{Z}}

\def	\R	{\mathbb{R}}
\def	\RR	{\mathbb{R}}
\def	\C	{\mathbb{C}}
\def	\CC	{\mathbb{C}}
\def	\CP	{\mathbb{C}\mathbb{P}}

\def    \ssminus{\smallsetminus}
\def	\half	{\frac{1}{2}}
\def	\del	{\partial}
\def	\eps	{\epsilon}
\def	\ol	    {\overline}
\def    \Fbar   {\ol{F}}
\def	\Hbar	    {{\ol{H}}}

\def	\bbar	    {{\ol{b}}}
\def	\xibar	    {{\ol{\xi}}}

\def    \tgamma     {{\tilde{\gamma}}}

\def    \PU         {\operatorname{PU}}

\def	\inv	{^{-1}}

\def \Sp {\text{Sp}}

\swapnumbers

\begin{document}

\title[Shortening Hofer length of circle actions]
{Shortening the Hofer length of Hamiltonian circle actions}

\author [Y.\ Karshon]{Yael Karshon}
%
%
\address{Department of Mathematics, University of Toronto,
40 St.George Street, Toronto ON M5S 2E4 Canada}
\email{karshon@math.toronto.edu}


\author [J.\ Slimowitz]{Jennifer Slimowitz}
\email{jslimowi@nsf.gov}

\thanks{Part of this research was conducted at the
Hebrew University of Jerusalem 
and was partially supported by the Israel Science Foundation,
founded by the Academy of Sciences and Humanities.
Jennifer Slimowitz acknowledges support by an NSF/NATO Postdoctoral 
Fellowship during the time that this work was done.}

\begin{abstract}
A Hamiltonian circle action on a compact symplectic manifold
is known to be a closed geodesic with respect to the Hofer metric
on the group of Hamiltonian diffeomorphisms.
If the momentum map attains its minimum or maximum 
at an isolated fixed point
with isotropy weights not all equal to plus or minus one,
then this closed geodesic can be deformed into a loop of
shorter Hofer length.  In this paper 
we give a lower bound for the possible amount 
of shortening, and we give a lower bound for the index
(``number of independent shortening directions"). 
If the minimum or maximum is attained along a submanifold $B$,
then we deform the circle action into a loop of shorter Hofer length
whenever the isotropy weights have sufficiently large absolute values
and the normal bundle of $B$ is sufficiently un-twisted.
\end{abstract}

\maketitle

\setcounter{tocdepth}{1}
\tableofcontents

\section{Introduction}
\labell{sec:intro}

A Hamiltonian circle action on a compact symplectic manifold $(M,\omega)$
can be viewed as a loop in the group $\Ham(M,\omega)$
of Hamiltonian diffeomorphisms, parametrized by $\R/\Z$.
We are concerned with deforming this loop,
through loops in $\Ham(M,\omega)$, into a loop whose Hofer length 
is smaller.  In short, we call such a deformation ``shortening".
The reader may go directly to Section~\ref{sec:new} 
for the precise statements of the new results that we prove in this paper.

\smallskip

One motivation for such a study comes from Riemannian geometry,
where one studies the action functional 
$\gamma \mapsto \int_0^1 \| \dot{\gamma}(t) \|^2 dt$ 
on spaces of loops or paths with fixed endpoints 
in a compact Riemannian manifold.
Critical points are geodesics; the index of a critical point 
is finite and is equal to the number (counting multiplicity) 
of conjugate points along the geodesic;
every homotopy class contains a minimal geodesic.
Morse theory in this context has played a major role in the study 
of the topology of compact Lie groups.  See Milnor's book~\cite{milnor}.

The infinite dimensional group $\Ham(M,\omega)$ is a central object 
in symplectic geometry.  It has a natural Finsler metric,
introduced by Hofer, which induces a non-trivial distance function.
The length of a path $\gamma \colon [0,1] \to \Ham(M,\omega)$
with respect to this Finsler metric is called the \emph{Hofer length}
of $\gamma$.  It is equal to the sum of the \emph{positive Hofer length}
and the \emph{negative Hofer length}; see section \ref{sec:preliminaries}.
%

There have been various attempts to extend to this setting
notions that had been applied successfully to the action functional
of Riemannian geometry.  The situation here is more difficult
and it is not yet clear to what extent the analogy can be carried out.
But the partial results that exist are beautiful and deep.
A main theme is to express Morse properties of Hofer's functional
on a space of paths in $\Ham(M,\omega)$
in terms of the dynamics of the flows on $M$ given by these paths.

There are several notions of ``geodesic path" in $\Ham(M,\omega)$.
With all of these notions, every time-independent Hamiltonian flow is geodesic.
(But not every geodesic is time-independent.)
%
%
Thus, every Hamiltonian circle action is a geodesic loop in $\Ham(M,\omega)$.

Starting with a Hamiltonian circle action 
$\{\phi_t \colon M \to M \}_{0 \leq t \leq 1}$
on a compact connected symplectic manifold,
here is a sample of known results.
Considering the opposite circle action gives the analogous results 
with ``maximum" replaced by ``minimum" and ``positive Hofer length"
replaced by ``negative Hofer length".

\begin{enumerate}
\item \labell{semifree}
Suppose that the circle action is semi-free near the set where the 
momentum map is maximal.  
(\emph{Semi-free} means that only the trivial group and the entire circle
occur as stabilizers.)
Then its positive Hofer length is minimal
among homotopic loops.
See~\cite{sav1}, the paragraph after Thm.~1.21,
whose proof relies on \cite[Thm.~1.9]{MT}.

In particular, if the circle action is semi-free on $M$ then it is
length-minimizing among homotopic loops in $\Ham(M,\omega)$.
This was already shown by McDuff-Slimowitz \cite[Theorem~1.4]{MSl}.

Another consequence is that,
if the circle action is semi-free near the set where the 
momentum map is maximal, then the circle action is non-contractible
in $\Ham(M,\omega)$.
This was already shown by McDuff-Tolman \cite[Theorem~1.1]{MT}.

\item \labell{Ust}
Suppose that the momentum map attains its maximum at an isolated 
fixed point and that the circle action is not semi-free near the maximum.
Then the circle action can be deformed to a loop in $\Ham(M,\omega)$
whose Hofer length is smaller. 
See Ustilovsky \cite[Theorem~1.2.E]{Us}.

\item \labell{trick}
If $n$ is a sufficiently large integer, 
then the non-effective circle action $\{\phi_{nt}\}_{0 \leq t \leq 1}$
can be deformed into a loop of smaller Hofer length.
See Polterovich~\cite[p.~86]{Pol}.

Let $B$ be the set where the momentum map attains its maximum.
If there exists a Hamiltonian diffeomorphism $g$ 
such that $g(B) \cap B = \emptyset$,
then the non-effective circle action $\{\phi_{2t}\}_{0 \leq t \leq 1}$ 
can be deformed to a loop in $\Ham(M,\omega)$ of smaller Hofer length.
See Polterovich \cite[Theorem~8.2.H]{Pol}.
\end{enumerate}

Many results in the literature apply to spaces of paths 
in $\Ham(M,\omega)$ with fixed endpoints
and to geodesic loops or paths that are possibly time-dependent.
We refer the reader to the papers
by Hofer~\cite{H}, Bialy-Polterovich~\cite{BP}, Yi Ming Long~\cite{long}, 
Siburg~\cite{Si}, Ustilovsky~\cite{Us}, 
Lalonde-McDuff~\cite{LM,LM:local},
McDuff-Slimowitz~\cite{MSl}, 
McDuff~\cite{McD:variants},
and Savelyev~\cite{sav1,sav2}.
Polterovich's book~\cite{Pol} contains systematic explanations
of many results in the subject.

\begin{Example}
Let $e^{2 \pi i t} \in S^1$ act on $\CP^2$ 
by $[u,z,w] \mapsto [u,e^{-2\pi i t}z,e^{2 \pi i t} w]$,
with the Fubini Study form normalized such that the momentum map is
%
%
$(|w|^2-|z|^2)/(|u|^2+|z|^2+|w|^2) + \text{constant}$.
Take an equivariant blow-up at the point $[0,\star,0]$
where the momentum map is minimal.
This yields a Hamiltonian circle action on a symplectic
manifold $(M,\omega)$ whose isotropy weights
where the momentum map is maximal are $-1,-2$ 
and whose isotropy weights 
where the momentum map is minimal are $1,1$. 
It is an effective circle action that represents a non-trivial
element of $\pi_1(M,\omega)$ (by the above item \eqref{semifree},
applied to the minimum)
but that can be deformed into a shorter loop 
in $\Ham(M,\omega)$ (by the above item \eqref{Ust}, applied
to the maximum).
\end{Example}

Typically, to show that a geodesic cannot be deformed
to a loop of smaller Hofer length requires ``hard" holomorphic curve
techniques,
whereas to show that a geodesic can be shortened
is possible with ``soft" techniques. 

In this paper we use ``soft" techniques to explore the behaviour 
of the Hofer length functional as we deform a Hamiltonian circle action 
within the space of loops in $\Ham(M,\omega)$. 
Starting with a circle action is a very restrictive assumption,
but it allows us to obtain new results that do not follow
from general existing results.  Our new results are sketched in the 
abstract and stated in Section~\ref{sec:new}.
In brief,
\begin{enumerate}
\item[$\star$]
Our Theorem~\ref{thm:isolated} is the first quantitative shortening
estimate for effective Hamiltonian circle actions.
\item[$\star$]
Our Theorem~\ref{thm:nonisolated}, Theorem~\ref{shorten w subbundle},
and Corollary~\ref{cor} are the first shortening results 
that apply when the extrema of the momentum map are not necessarily 
displaceable and the circle action is effective.  
\item[$\star$]
Our Theorem~\ref{thm:savelyev} gives a lower bound on the index 
of the Hofer length functional. It confirms a prediction 
of Yasha Savelyev. (But see Remark~\ref{rk:foot}.)
\end{enumerate}

\bigskip

We now give an overview of the paper.

\smallskip

In Section \ref{sec:preliminaries} we review standard background material.
The purpose of this review is to remind
equivariant symplectic geometers of facts about time-dependent
Hamiltonian flows and to remind symplectic topologists
of facts about momentum maps.
In Section \ref{sec:new} we state our results, give examples,
and pose further questions.

\smallskip

In Sections~\ref{sec:polterov}--\ref{sec:disjoin family}
we develop the main tools
for the later proofs of our main results.
We warm up by describing, in Section \ref{sec:polterov},
Polterovich's shortening trick for non-effective circle actions
(item~\eqref{trick} in the above sample of known results).
In our later proofs we use variations of this Polterovich trick;
we give such variations in Section \ref{sec:polterov2}:
take a Hamiltonian circle action with momentum map $H$,
and assume that, near the set where $H$ is maximal, we can write
$H = K+F$ where $K$ and $F$ also generate circle actions.
Also assume that the sets where $K$ and $F$ are maximal have non-empty
intersection, and that there exists a symplectic isotopy
that disjoins the first of these sets from the second.
Finally, assume that the infimum of the sum of $K$ and $F$ is equal
to the sum of their infima.
We deduce that the circle action can be deformed
to a shorter loop in $\Ham(M,\omega)$,
and we estimate the amount by which this loop is shorter.
Our actual statement is slightly more technical; 
see Lemma~\ref{lem polterov}.
The qualitative part of this result can be strengthened if $H$ decomposes, 
near its maximum, into a sum of more than two momentum maps;
see Lemma~\ref{lem polterov k}.
Sections \ref{sec:disjoin disc} and \ref{sec:disjoin family}
contain variations on the simple fact that a disc of area $A_1$
can be disjoined from a disc of area $A_2$
inside a disc of area greater than $A_1 + A_2$.
See Lemmas~\ref{lem:btA} and~\ref{CtimesW}.

\smallskip

In Section~\ref{sec:isolated} we prove our first new result,
Theorem~\ref{thm:isolated}:
when the momentum map attains its maximum at an isolated fixed point
and one of the isotropy weights has absolute value $\geq 2$,
we give a lower bound for the amount of possible shortening,
in terms of the gap between the two largest critical values
of the momentum map.
The proof consists of the following steps.
There exists an equivariant Darboux chart whose image
is the entire set of points where the momentum map
is larger than the second-largest critical value;
this set can be identified with an ellipsoid in $\C^n$.
We view this ellipsoid as a family of discs in $\C$
parametrized by points in a subset of $\C^{n-1}$,
where the circle action on the $\C$ component is non-effective.
From the results of Section \ref{sec:polterov2} 
we get shortening by an amount that 
depends on
the size of a family of discs 
that can be disjoined from another such family
through a Hamiltonian isotopy supported in the ellipsoid.   
Such an estimate, in turn, is obtained 
from the results of Section~\ref{sec:disjoin family}.
For effective circle actions,
these are the first shortening results that are quantitative.

\smallskip

In Sections~\ref{sec:nbhd of B} and~\ref{sec:nonisolated}
we prove our second new result (more precisely, set of results).
We now allow the momentum map to attain its maximum 
along a submanifold $B$ of positive dimension. 
To describe a neighbourhood of $B$ in $M$,
we use Sternberg's minimal coupling procedure
and Weinstein's symplectic tubular neighbourhood theorem.
Let $E$ denote the normal bundle of $B$ in $M$.
In the simplest case, we may assume that $E$ has $k$ sections
that are nowhere all vanishing, where $k+1$ is the smallest
absolute value of an isotropy weight in $E$.
In the general case, we make a similar assumption
for a \emph{sub-bundle} $E'$ of $E$.
The idea of the proof is to express the circle action
as a composition of $k+1$ circle actions
on a neighbourhood of $B$ and to then apply the results
of Section~\ref{sec:polterov2}.
For the precise results, and for some special cases, 
see Theorems~\ref{thm:nonisolated} and~\ref{shorten w subbundle}
and Corollary~\ref{cor}.
These are the first shortening results when the maximum
is attained along a manifold that is not necessarily displaceable
and the circle action is effective.

\smallskip

Finally, in Section~\ref{sec:savelyev}, we prove our third new result,
Theorem~\ref{thm:savelyev}:
we give a lower bound for the index of the Hofer length functional,
assuming that the momentum map attains its maximum at an isolated fixed point.
Composing with an equivariant Darboux chart, we may work with 
a linear circle action on $\C^n$
that rotates the coordinates of $\C^n$ with speeds $-k_1,\ldots,-k_n$,
where the $k_j$s are positive integers.
We write the action on the $j$th coordinate
as a composition of $k_j$ terms, $\phi_t \circ \ldots \circ \phi_t$,
where $\phi_{-t}$ generates the scalar multiplication circle action.
We obtain a family of shortenings by applying Polterovich's trick
of Section~\ref{sec:polterov} to the different factors
of these compositions: in each coordinate, we conjugate all but 
one of the $\phi_t$s by a translation of~$\C$.  Varying over the
two dimensional family of translations of $\C$,
independently in each one of the $\sum(k_j-1)$ conjugations,
we obtain a family of loops in $\Ham(M,\omega)$
that depends on $2 \sum (k_j-1)$ real parameters
and such that the original circle action
corresponds to the origin in the parameter space $\R^{2 \sum (k_j-1)}$.
An explicit computation then shows that the
Hofer length of these loops, as a function on $\R^{2 \sum (k_j-1)}$,
achieves its maximum at the origin, and its Hessian is negative definite
at the origin.
This confirms one direction of a prediction for the index
that was made by Yasha Savelyev.

\begin{Remark} \labell{rk:foot}
While completing this paper, we learned that Savelyev \cite{sav3} 
can now confirm his conjecture with the methods of Ustilovsky \cite{Us}.
\end{Remark}

\bigskip\noindent
This paper is the outgrowth of a joint project that took place in 
the years 1999-2001.  Yael Karshon takes the responsibility for any errors 
that might be present.

\bigskip\noindent\textbf{Acknowledgement}
We thank Leonid Polterovich for inspiring conversations
and helpful suggestions.  We thank Yasha Savelyev for
suggesting to compute the index of the Hofer length functional
and telling us his conjecture for its value.
We also thank Fabian Ziltener and Brian Lee for helpful discussions.
Finally, we are grateful to the referee for his or her extremely 
thorough report and numerous helpful comments.

\section{Preliminaries}
\labell{sec:preliminaries}

The purpose of this section is to remind equivariant symplectic geometers
of facts about time-dependent Hamiltonian flows and to remind 
symplectic topologists of facts about momentum maps.

\subsection*{Hamiltonian isotopies.}
Let $(M,\omega)$ be a symplectic manifold.
A \emph{smooth path} of symplectomorphisms
is a map $t \mapsto \psi_t$, from an interval~$I$
to the group of symplectomorphisms of $(M,\omega)$,
such that the map $(t,m) \mapsto \psi_t(m)$ 
from $I \times M$ to $M$ is smooth.
Such a path is a \emph{Hamiltonian isotopy}
if its velocity vector field, defined by
$\frac{d}{dt} \psi_t = X_t \circ \psi_t $,
is generated by a smooth function
$H \colon I \times M \to \R$,
by Hamilton's equations
$ d H_t = - \iota(X_t) \omega$,
where $H_t(\cdot) = H(t,\cdot)$.
The function $H = \{ H_t \}_{t \in I}$
is called \emph{the Hamiltonian}.
If the Hamiltonian isotopy $\{ \psi_t \}$ is compactly supported, 
its \emph{Hofer length} is
$$ \length \{ \psi_t\}_{t \in I} =
\int\limits_I \left( \max_{x \in M} H_t(x) - \min_{x \in M} H_t(x) 
	 \right) dt.$$
A symplectomorphism is \emph{Hamiltonian} if it can be connected
to the identity map by a Hamiltonian isotopy.

\subsection*{The positive and negative Hofer length functionals}
Suppose that $M$ is compact.
The Hamiltonian $\{ H_t \}$ is \emph{normalized} if 
$$\int_M H_t \, \omega^n = 0 \quad \text{ for all } t.$$
The \emph{positive Hofer length functional} associated to a 
Hamiltonian isotopy
$\{ \varphi_t \}_{a \leq t \leq b}$ in $\Ham(M,\omega)$
is the number
$$ \ell_+(\varphi) = \int_a^b \max\limits_M H_t \, dt $$
where $\{H_t\}_{a \leq t \leq b}$ is the normalized Hamiltonian
that generates the isotopy.
The \emph{negative Hofer length functional} is, similarly,
$$ \ell_-(\varphi) = \int_a^b \max\limits_M -H_t \, dt.$$
The Hofer length of $\varphi$ is then 
$\ell_+(\varphi) + \ell_-(\varphi)$.

\subsection*{Conjugation}
If $\{ \psi_t \}$ is a Hamiltonian isotopy,
generated by $\{ H_t \}$,
and $b$ is a symplectomorphism,
then $ \{ b \psi_t b\inv \}$ is a Hamiltonian
isotopy, generated by $\{ H_t \circ b\inv \}$.
If $H_t$ is normalized then $H_t \circ b\inv$
is also normalized.

\subsection*{Reparametrization.}
Let $\{ \psi_t \}$ be a Hamiltonian isotopy, generated
by a function $\{ H_t \}$.
Let $t = t(\tau)$ be a smooth function.
Then $\psibar_\tau := \psi_{t(\tau)}$
is a Hamiltonian isotopy, generated by 
the product
$\Hbar_\tau := H_{t(\tau)} \frac{dt}{d\tau} $.

Indeed, let $\xi_t$ be the vector field 
that generates the isotopy, so that
$\frac{d}{dt} \psi_t = \xi_t \circ \psi_t$
and $dH_t = - \iota(\xi_t) \omega$.
Let $\xibar_\tau = \xi_{t(\tau)} \frac{dt}{d\tau} $.
Then 
$ d\Hbar_\tau
= dH_{t(\tau)} \frac{dt}{d\tau} 
= - \iota(\xibar_\tau) \omega $
and $\frac{d}{d\tau} \psibar_\tau
= \frac{d}{dt} \left. \psi_t \right|_{t=t(\tau)} \frac{dt}{d\tau} 
= \xibar_\tau \circ \psibar_\tau$.

\subsection*{Composition and inverse}
If $\{ \psi^K_t \}_{t \in I}$ is a Hamiltonian isotopy
generated by $\{ K_t \}_{t \in I}$,
and $\{ \psi^F_t \}_{t \in I}$ is a Hamiltonian isotopy
generated by $\{ F_t \}_{t \in I}$,
then $\{ \psi^K_t \circ \psi^F_t \}_{t \in I}$ is a Hamiltonian 
isotopy generated by $\{ K_t + F_t (\psi^K_t)\inv \}_{t \in I}$,
and $(\Psi^K_t)\inv$ is a Hamiltonian isotopy
generated by $-K_t \circ \Psi^K_t$.
If $K_t$ and $F_t$ are normalized, 
so are $K_t + F_t \left( \psi^K_t \right)\inv$
and $-K_t \circ \Psi^K_t$.

The group of Hamiltonian symplectomorphisms is denoted $\Ham(M,\omega)$.

\subsection*{Paths in the Hamiltonian group}
If $\{ \psi_t \}$ is a smooth path in the group
of symplectomorphisms of $(M,\omega)$,
and for each $t$ the symplectomorphism $\psi_t$ 
is in the subgroup $\Ham(M,\omega)$,
then $\{ \psi_t \}$ is a Hamiltonian isotopy.
This is a result due to Banyaga~\cite[p.190, prop. II.3.3]{Ba}.
In textbooks, the proof is often omitted or is intertwined 
with proofs of more difficult facts. We recall the proof.

Let $\xi_t$ be the velocity vector field,
defined by $\frac{d}{dt}\psi_t = \xi_t \circ \psi_t$.
The path $\{ \psi_t \}$ is a Hamiltonian isotopy
if and only if $\iota(\xi_t)\omega$ is exact for all $t$,
if and only if the $H^1(M;\R)$-valued integral
$\Flux(\{\psi_\tau\}_{0 \leq \tau \leq t})
:= \int_0^t [\iota(\xi_\tau)\omega] d\tau$
is zero for all $t$.
(We then set $H_t(x) = - \int_{x_0}^x \iota(\xi_t) \omega$.)

The evaluation of the flux on a smooth loop
$\gamma \colon S^1 \to M$ is equal to the 
symplectic area of the surface
$\tgamma \colon S^1 \times [0,t] \to M$
given by $\tgamma(s,\tau) := \psi_\tau(\gamma(s))$;
this follows from Stokes's formula.
It follows that the flux of a \emph{closed} loop 
of symplectomorphisms is a class in $H^1(M;\R)$
whose evaluation on elements of $H_1(M;\Z)$
takes values in the countable set
$ \left< \omega , H_2(M;\Z) \right> $.
Such classes make up a countable subgroup of $H^1(M;\R)$.

Suppose that $\{ \psi_t \}$ is a path of symplectomorphisms
and that each $\psi_t$ is Hamiltonian.
Then for each $t$ the path $\{ \psi_\tau \}_{0 \leq \tau \leq t}$
can be completed to a closed loop of symplectomorphisms
by a Hamiltonian isotopy.  This completion does not effect the flux.
So $t \mapsto \Flux(\{ \psi_\tau \}_{0 \leq \tau \leq t})$
is a continuous map that takes values in a countable subgroup
of $H^1(M;\R)$.
Such a map must be constant.
Since it is zero for $t=0$, it is zero for all $t$.

\subsection*{Hamiltonian circle actions}

Let $(M,\omega)$ be a symplectic manifold 
with a circle action generated by a Hamiltonian
$H \colon M \to \R$.
The \emph{momentum map} for the circle action
is the Hamiltonian $H$.  

The critical points of $H$
are the fixed points for the circle action.
Let $p \in M$ be such a point.
The linearized isotropy action on $T_pM$
is linearly equivariantly symplectomorphic
to $\C^n$ with the standard symplectic form
and with the circle acting as a subgroup of $(S^1)^n$,
namely, $e^{2 \pi i t} \cdot (z_1,\ldots,z_n)
= (e^{2 \pi i k_1 t} z_1 , \ldots , e^{2 \pi i k_n t} z_n )$. 
The integers $k_1,\ldots,k_n$, which measure the speeds at which
the linearized $S^1$-action rotates its eigenspaces,
are called the \emph{isotropy weights} at $p$.
They are unique up to permutation.
This isomorphism $T_pM \cong \C^n$ carries the
Hessian of $H$ to the quadratic form
$ \pi k_1 |z_1|^2 + \ldots + \pi k_n |z_n|^2 $.
So if $H$ attains its maximum at $p$ then all the isotropy weights 
at $p$ are non-positive.

The \emph{equivariant Darboux} theorem asserts that a neighbourhood
of $p$ in $M$ is equivariantly symplectomorphic
to a neighbourhood of the origin in $T_pM$.
In particular, if $p$ is an isolated fixed point
then the isotropy weights at $p$ are all non-zero.
Thus, at least one isotropy weight at $p$ has absolute value 
greater than one
if and only if every neighbourhood of $p$ contains a point
whose stabilizer is finite and non-trivial. 
Also, the set $\Bmax$ where $H$ is maximal
is a symplectic submanifold, and the circle acts on the 
fibres of its normal bundle with negative weights.

The momentum map $H$ is a Bott-Morse function with even indices
and coindices.
Assuming that $M$ is compact and connected, this implies that 
the level sets of $H$ are connected and that preimages of intervals 
are connected.  
Thus, the maximum of $H$ is attained along a connected component
of the fixed point set.
If a fixed point $p$ is a local maximum for $H$
then it is also a global maximum.
If a neighbourhood of $p$ in $M$
is equivariantly symplectomorphic to the ellipsoid
$\{ z \ | \ \sum k_j \pi |z_j|^2 < \eps \}$ 
with the momentum map $- \sum k_j \pi |z_j|^2$,
then this neighbourhood contains \emph{all}
the points in $M$ where $H$ is $\eps$-close to its maximum.

See~\cite{At,GS:convexity}.

\section{New results and further questions}
\labell{sec:new}

In this section we state our main results.
We prove them in Sections~\ref{sec:isolated},
\ref{sec:nonisolated}, and~\ref{sec:savelyev}.
We also give some examples, corollaries, and further questions. 

In each of these results, we start with a Hamiltonian circle action,
viewed as a loop in the group of Hamiltonian diffeomorphisms,
and we deform it to a loop whose Hofer length is smaller.
Moreover, throughout the deformation, the negative Hofer length
remains constant, and the positive Hofer length eventually becomes smaller
without ever becoming larger than the initial one.

\subsection*{Shortening on a manifold with an isolated maximum.}

\begin{Theorem} \labell{thm:isolated}
Let $(M,\omega)$ be a compact connected symplectic manifold
with a Hamiltonian circle action.
Suppose that the momentum map attains its maximum 
at an isolated fixed point
and that at least one of the isotropy weights at that point
has absolute value greater than one.
Let~$d$ be a positive number that is smaller than the gap
between the two largest critical values of the momentum map.
Then the circle action can be deformed to a loop in $\Ham(M,\omega)$ 
whose Hofer length is smaller than that of the original circle action
by $2d/9$.
\end{Theorem}

\begin{Remark} \labell{rk:isolated}
In Theorem~\ref{thm:isolated}, we can choose the deformation
such that the positive Hofer length never exceeds the initial one
and the negative Hofer length remains constant.
Specifically, let $H$ denote the momentum map, let $\pmax$ denote
the point where $H$ attains its maximum, let $\alpha$ denote the
gap between the two largest critical values of $H$,
and consider the subset of $M$ given by
\begin{equation} \labell{alpha nbhd'}
\{ m \ | \ H(m) > H(\pmax) - \alpha \}.
\end{equation}
We can choose the deformation such that 
the normalized deformed Hamiltonians remain greater than $H(\pmax)-\alpha$
on the set~\eqref{alpha nbhd'} and coincide with~$H$ outside this set.
\end{Remark}

\begin{Remark} \labell{rk:isolated'}
In Theorem~\ref{thm:isolated}, we actually obtain a better estimate 
than $\frac{2}{9}d$.
If one of the isotropy weights is even then we can shorten by $\frac{1}{4}d$.
Otherwise, let $-k_1$ be the weight whose absolute value is largest;
then we can shorten by $\frac{1}{4} ( 1 - \frac{1}{k_1^2}) d$.
When $k_1=3$, this becomes $\frac{2}{9}d$.
\end{Remark}

This is the first quantitative result for shortening 
of loops in $\Ham(M,\omega)$ 
that are effective circle actions.
We prove it in Section~\ref{sec:isolated}.

The estimate of Theorem~\ref{thm:isolated} is far from sharp,
as we see in the following example.

\begin{Example} \labell{not optimal}
Let $a,b$ be positive integers that are relatively prime.
Let $e^{2\pi i t} \in S^1$ act on $\CP^2$ by 
$[u,z,w] \mapsto [u, e^{- 2 \pi i t a} z , e^{2 \pi i t b} w]$,
with the Fubini-Study form normalized such that 
the momentum map is
$ (-a|z|^2+b|w|^2) / (|u|^2+|z|^2+|w|^2) + \text{constant}$.
The Hofer length of this action is $a+b$.
The isotropy weights at the maximum are $-b,-(a+b)$;
the isotropy weights at the minimum are $a,a+b$.
Applying the shortening result of Theorem~\ref{thm:isolated}
and the analogous result with maximum replaced by minimum,
for any $\eps > 0$
we can deform the circle action to a loop of Hofer length
$< 7/9 (a+b) + \eps$. 
This upper bound is greater than one.
But the action can be deformed to a loop of Hofer length one or zero,
so the bound is not optimal.
(To see that the action can be deformed to a loop 
of Hofer length one or zero, 
notice that the action extends to an action of $\PU(3)$,
that $\pi_1(\PU(3)) \cong \Z/3\Z$,
and that the non-trivial elements of $\pi_1(\PU(3))$ 
are represented by the actions
$[u,z,w] \mapsto [u, e^{-2 \pi i t} z , w]$
and
$[u,z,w] \mapsto [u, z , e^{2 \pi i t} w]$,
which have Hofer length equal to one.)
\end{Example}

\begin{Remark} \labell{on not optimal}
The above example is inspired by Proposition~1.3 of the paper~\cite{MT2}
of McDuff and Tolman:  
if $G$ is a compact Lie group with trivial centre
and $M$ is a coadjoint orbit of $G$,
then every non-trivial element of $\pi_1(G)$
is represented by a sub-circle of $G$ that acts on $M$ semi-freely.
\end{Remark}

\subsection*{Shortening on a manifold with an arbitrary maximum.}\ \\

Recall that, on a compact
connected symplectic manifold with a Hamiltonian circle action,
the set where the momentum map is maximal is a connected symplectic
submanifold.

In the following theorem, a subset $B$ of a symplectic manifold $M$
is \emph{symplectically $k$-displaceable} in a neighbourhood $V$
if there exist $k$ symplectomorphisms $b_1,\ldots,b_k$ of $M$,
each connected to the identity through a path of symplectomorphisms 
supported in $V$, such that
$B \cap b_1(B) \cap \ldots \cap b_k(B) = \emptyset$.

In the symplectic literature, ``$B$ is displaceable" often means
that there exists a Hamiltonian symplectomorphism $g$
such that $B \cap g(B) = \emptyset$.
Thus, a set can be symplectically $k$-displaceable in every neighbourhood
without being displaceable (even when $k = 1$).

{
\begin{Theorem} \labell{thm:nonisolated}
Let $(M,\omega)$ be a compact connected symplectic manifold
with a Hamiltonian circle action.
Let $\Bmax$ be the set where the momentum map is maximal.
Let $-k_1,\ldots,-k_s$ denote the distinct weights
for the circle action on the normal bundle of $\Bmax$,
and let $k := \min \{ k_1,\ldots,k_s \} - 1$.
Suppose that $\Bmax$ is symplectically $k$-displaceable
in every neighbourhood.
Then the circle action can be deformed through loops in $\Ham(M,\omega)$
into a loop of smaller Hofer length.
\end{Theorem}

\begin{Remark} \labell{rk:nonisolated}
In Theorem~\ref{thm:nonisolated},
we can choose the deformation such that the positive Hofer length 
never exceeds the initial one and the negative Hofer length remains constant.
Specifically, let $H$ denote the momentum map. 
For every positive number $\alpha$,
we can choose the deformation such that the normalized deformed Hamiltonians
remain between $H(\Bmax)-\alpha$ and $H(\Bmax)$ on the set
$\{ m \in M \ | \ H(m) > H(\Bmax) - \alpha \}$
and coincide with $H$ outside this set.
\end{Remark}
}

In practice, the easiest way to show that a symplectic submanifold
is symplectically $k$-displaceable in every neighbourhood is to show 
that its normal bundle has $k$ sections that are nowhere all vanishing.
($k$-displaceability then follows 
from the symplectic tubular neighbourhood theorem~\eqref{lnf}
and Lemma~\ref{vf}.)
More generally, we have the following theorem.

\begin{Theorem}  \labell{shorten w subbundle}
Let $(M,\omega)$ be a compact connected symplectic manifold
with a Hamiltonian circle action.
Let $\Bmax$ be the set where the momentum map is maximal.
Let $E' \to \Bmax$ be an $S^1$-invariant subbundle
of the normal bundle to $\Bmax$ in $M$,
and let $-k'_1, \ldots, -k'_{s'}$ be the distinct weights
for the circle action on $E'$.
Let 
$$ k' = \min \{ k'_1, \ldots, k'_{s'} \} -1 . $$
Suppose that 
\begin{equation}  \labell{condition on subbundle}
\text{$E'$ has $k'$ sections that are nowhere all vanishing.}
\end{equation}
Then the circle action can be deformed through loops in $\Ham(M,\omega)$
into a loop of smaller Hofer length.
\end{Theorem}

\begin{Remark} \labell{rk:shorten w subbundle}
In Theorem~\ref{shorten w subbundle},
we can choose the deformation such that the positive Hofer length 
never exceeds the initial one and the negative Hofer length remains constant.
Specifically, let $H$ denote the momentum map. 
For every positive number $\alpha$,
we can choose the deformation such that the normalized deformed Hamiltonians
remain between $H(\Bmax)-\alpha$ and $H(\Bmax)$ on the set
$\{ m \in M \ | \ H(m) > H(\Bmax) - \alpha \}$
and coincide with $H$ outside this set.
\end{Remark}

Theorem~\ref{shorten w subbundle} has the following corollaries.

\begin{Corollary} \labell{cor}
Let $(M,\omega)$ be a compact symplectic manifold
with a Hamiltonian circle action.
Let $\Bmax$ be the subset of $M$ where the momentum map
attains its maximum.
Let $E$ denote the normal
bundle of $\Bmax$ in $M$ and $-k_1,\ldots,-k_s$ the distinct weights
for the circle action on $E$.  
Suppose that one of the following conditions holds.

\begin{enumerate}

\item 
The normal bundle of $\Bmax$ in $M$
has a trivial sub-bundle on which the circle acts with a weight 
of absolute value $\geq 2$.

(In particular, this holds under the assumptions of 
Theorem~\ref{thm:isolated}.)

\item \labell{prevcase}
We have
$$k_j > 1 + \frac{\dim \Bmax}{\codim \Bmax}$$
for all $j \in \{ 1,\ldots,s \}$.

\item 
After a possible relabeling, assume that  $k_1 < \ldots < k_s$.
Let $E_j$ denote the sub-bundle where the circle acts with weight $-k_j$.
There exists $j \in \{ 1 , \ldots , s \}$
such that 
$$ (k_j-1) \rank(E_j \oplus \ldots \oplus E_s) > \dim \Bmax .$$ 
(For $j=1$, this condition amounts to case~\eqref{prevcase}.)
\end{enumerate}
Then the circle action can be deformed through loops in $\Ham(M,\omega)$
to a loop of smaller Hofer length.
\end{Corollary}

Moreover, we can choose a deformation with the properties
described in Remark~\ref{rk:shorten w subbundle}.

\begin{proof}   [Proof of Corollary~\ref{cor}]
In each of the three cases, \eqref{condition on subbundle} holds
for an appropriate choice of sub-bundle $E'$:
\begin{enumerate}
\item
Let $E'$ be a trivial sub-bundle of $E$
on which $S^1$ acts with a weight $-k'_1$ 
of absolute value $k_1 \geq 2$. 
Then $k' = k'_1 - 1 \geq 1$,
and~\eqref{condition on subbundle} holds because $E'$ is trivial.
\item
Let $E'=E$. Then $\rank E' = \codim \Bmax$.
Our assumption on the weights give
$k' > \frac{\dim \Bmax}{\codim \Bmax}$.
So $k' \rank E' > \dim \Bmax$,
and~\eqref{condition on subbundle} holds for generic sections.
\item
Let $E' = E_j \oplus \ldots \oplus E_s$.
Then $k' = k_j - 1$.  
The assumption on the weights gives
$k' \rank E' > \dim \Bmax$,
and again~\eqref{condition on subbundle} holds for generic sections.
\end{enumerate}
In each of these cases, the shortening result
follows from Theorem~\ref{shorten w subbundle}.
\end{proof}

\begin{Remark}
Let $k$ and $k'$ be as in Theorems~\ref{thm:nonisolated}
and~\ref{shorten w subbundle}.
Existence of $k'$ sections of $E'$ that are nowhere all vanishing
does not necessarily imply the existence of $k$ sections of $E$
that are nowhere all vanishing, because $k'$ can be larger than $k$.
So Theorem~\ref{shorten w subbundle} does not follow 
from Theorem~\ref{thm:nonisolated}.
\end{Remark}

\begin{Remark}
Informally, the assumptions of Theorems~\ref{thm:nonisolated} 
and~\ref{shorten w subbundle} and Corollary~\ref{cor}
mean that the normal bundle is ``sufficiently untwisted".
\end{Remark}

\begin{Remark} \labell{visible}
%
%
There seems to be an interesting relationship 
with the work of McDuff-Tolman.
In the terminology of Section~1.1 of McDuff-Tolman's paper~\cite{MT}, 
the component $\Bmax$ is \emph{homologically visible} if and only if 
the Euler class of the \emph{obstruction bundle}
$(E_1 \otimes \C^{k_1-1}) \oplus \ldots \oplus
(E_s \otimes \C^{k_s-1})$ is non-zero.
Here, $k_1,\ldots,k_s$ are the weights of the circle action
on the normal bundle, and $E_j$ is the subbundle on which
the circle acts with weight $k_j$.
By Theorem~1.2 of~\cite{MT}, if $\Bmax$ is homologically visible
then the circle action is non-contractible in $\Ham(M,\omega)$.
(We are not sure whether $\Bmax$ being homologically visible
further implies that the circle action is minimal among homotopic loops.
Proposition~3.3 of~\cite{MT} might be relevant.)
%
%
The condition~\eqref{condition on subbundle} 
of our Theorem~\ref{shorten w subbundle}
implies that $\Bmax$ is not homologically visible.
($k'$ sections of $E'$ that are nowhere all vanishing
can be considered as a non-vanishing section of $E' \otimes \C^{k'}$.
Writing $E' = E'_1 \oplus \ldots \oplus E'_s$ where the circle
acts on $E'_j$ with weight $k_j$, we get a nonvanishing section
of $(E'_1 \otimes \C^{k'}) \oplus \ldots \oplus (E'_{s'} \otimes \C^{k'})$.
Because $k' \leq k'_j - 1$ for all $j=1,\ldots,s'$,
we get a non-vanishing section of
$(E'_1 \otimes \C^{k'_1-1}) \oplus \ldots 
\oplus (E'_{s'} \otimes \C^{k'_s-1})$.
Because this is a subbundle of the obstruction bundle,
the obstruction bundle has a non-vanishing section, 
and hence its Euler class is zero.)
%
%
\end{Remark}

\smallskip

We prove Theorems~\ref{thm:nonisolated} and~\ref{shorten w subbundle}
in Section~\ref{sec:nonisolated}.  
These are the first shortening results that apply to situations
when the maximum is not necessarily displaceable. 

\begin{Example}
Let $e^{2 \pi i t} \in S^1$ act on $\CP^4$ by
$$[z_0,z_1,z_2,z_3,z_4] \mapsto
[z_0,z_1, e^{-2 \pi i t} z_2, e^{-2 \pi i t \cdot 2} z_3,
e^{-2 \pi i t \cdot 3} z_4],$$
with momentum map
$( -|z_2|^2 - 2|z_3|^2 -3|z_4|^2 ) / 
( |z_0|^2 + |z_1|^2 + |z_2|^2 + |z_3|^2 + |z_4|^2 ) + \text{constant}.$
The momentum map attains its maximum along the submanifold 
$\Bmax = \{ [z_0,z_1,0,0,0] \}$, which is isomorphic to $\CP^1$.
The negative isotropy weights along $\Bmax$ are $-1,-2,-3$.
The momentum map attains its minimum at the isolated fixed point 
$[0,0,0,0,z_4]$.
The isotropy weights at the minimum are $1,2,3,3$.
In the notation of item~(3) of Corollary~\ref{cor}, applied at $\Bmax$,
we have the following values.
\begin{center}
\begin{tabular}{c|c|c}
$j$ & $k_j$ & $\rank E_j$ \\
\hline  
1 & 1 & 2 \\
2 & 2 & 2 \\
3 & 3 & 2 
\end{tabular}
\end{center}
We can apply item~(3) of Corollary~\ref{cor} to either
$E_2 \oplus E_3$ or $E_3$.  
Also applying Theorem~\ref{thm:isolated} at the minimum,
we conclude that 
the circle action can be deformed to a loop 
whose positive Hofer length and negative Hofer length are both smaller.
Nevertheless, the circle action is not contractible in $\Ham(M,\omega)$:
it is homotopic (in $\PU(5)$) to the circle action
$[z_0,z_1,z_2,z_3,z_4] \mapsto [e^{-2\pi i t} z_0, z_1,z_2,z_3,z_4,z_5]$,
which is semi-free and thus non-contractible
(by item~\eqref{semifree} on p.~\pageref{semifree}).
\end{Example}

\subsection*{The index of the Hofer positive length functional.}

\begin{Theorem} \labell{thm:savelyev}
Let $(M,\omega)$ be a $2n$ dimensional
compact connected symplectic manifold with a Hamiltonian
circle action.  
Suppose that the momentum map attains its maximum at an isolated fixed point.
Let $-k_1,\ldots,-k_n$ be the isotropy weights at the maximum,
with possible repetitions.
Then there exists a neighbourhood $D$ of the origin 
in $\R^{\sum 2(k_j - 1)}$, and, for each $\lambda \in D$, a loop
$\{ \psi^{(\lambda)}_t \}_{0 \leq t \leq 1}$ in $\Ham(M,\omega)$,
such that the following properties hold.
For $\lambda=0$, the loop $\{\psi^{(0)}_t\}_{0 \leq t \leq 1}$ 
is the given circle action.
The function $\lambda \mapsto \length( \{ \psi^{(\lambda)}_t \} )$ is smooth, 
$\lambda=0$ is a critical point of this function,
and the Hessian of this function at $\lambda=0$ is negative definite.
\end{Theorem}

\begin{Remark} \labell{rk:savelyev}
In fact, we can choose the deformation such that the negative 
Hofer length remains constant.  Specifically, let $H$ denote
the momentum map, and let $\pmax$ denote the point where $H$ attains
its maximum. For every positive number $\alpha$,
we can choose the deformation such that the normalized deformed Hamiltonian
remains between $H(\pmax)-\alpha$ and $H(\pmax)$ on the set
$\{ m \in M \ | \ H(m) > H(\pmax)-\alpha \}$
and coincides with $H$ outside this set.
\end{Remark}

Savelyev \cite{sav1,sav2}, in the setup of Theorem~\ref{thm:savelyev},
predicted that the index of the positive Hofer length functional 
on the space of loops in $\Ham(M,\omega)$
at the given Hamiltonian circle action is equal to $\sum 2(k_i-1)$.
Theorem~\ref{thm:savelyev} is a sense in which the index 
of the positive Hofer length functional
is \emph{at least} $\sum_i 2(k_i - 1)$.
So Theorem~\ref{thm:savelyev} confirms the $\geq$ direction 
of Savelyev's prediction.
(But see Remark~\ref{rk:foot}.)

\subsection*{Further questions}
\begin{enumerate}
\item
Are there examples where the estimate of Theorem~\ref{thm:isolated} 
is sharp?
\item
Can the estimate of Theorem~\ref{thm:isolated} be improved
in the presence of high isotropy weights?

As noted by Polterovich,
this question is also interesting for non-effective
circle actions obtained by iterating an effective circle action,
as it may provide information on the asymptotic norm 
of the homotopy class of the effective circle action (see~\cite{Pol}).

\item
Under the assumptions of Theorem~\ref{thm:nonisolated}
or of Theorem~\ref{shorten w subbundle},
can one obtain a quantitative estimate?
\item
In Theorem~\ref{shorten w subbundle},
can the condition~\eqref{condition on subbundle}
be replaced by the condition that $\Bmax$ 
not be homologically visible?
(Cf.~Remark \ref{visible}.)
\item
Can Theorem~\ref{thm:savelyev} be extended to cases in which 
the maximum is attained along a manifold of positive dimension?

\item
Savelyev's prediction for the index of the positive Hofer length functional
is an even number. Polterovich expects the index to be even
in additional cases.
On a compact connected manifold, a Morse function with even indices 
and coindices has nice topological properties:
every local minimum is a global minimum,
and the level sets are connected.
Palais-Smale-Morse theory does not apply to the space
of loops in $\Ham(M,\omega)$, but we may seek analogies.  

If a loop in the Hamiltonian group is a local minimum for the
positive Hofer length functional, 
is it also a global minimum in its homotopy class?

Suppose that a homotopy class in $\pi_1(\Ham(M,\omega))$ contains two
circle actions that are semi-free near the set where the momentum
map is maximal.   (Cf.\ item~(1) on p.~\pageref{semifree}.)
Can we deform one to the other through loops of constant
positive Hofer length?
\end{enumerate}

\section{Shortening a non-effective circle action}
\labell{sec:polterov}

To warm up, we describe a procedure,
which we learned from Leonid Polterovich,
for shortening a non-effective Hamiltonian circle action
by displacing the set where its momentum map is maximal.

Let $H \colon M \to \R$ be a time independent Hamiltonian.
Suppose that $\frac{1}{2} H$ generates a circle action, 
$\{\phi_t\}_{0 \leq t \leq 1}$, so that $H$ generates the
non-effective circle action
\begin{equation} \labell{phikt}
\psi^H_t = \phi_{2t} = \phi_t \circ \phi_{t}. 
\end{equation}

Let $b \colon M \to M$ be a symplectomorphism. 
Then
\begin{equation} \labell{shorter loop}
\psi_t^\Hbar := \{ \phi_t \, b \, \phi_{t} \, b\inv \}_{0 \leq t \leq 1}
\end{equation}
is a loop in $\Ham(M)$, based at the identity,
and generated by the Hamiltonian 
\begin{equation} \labell{Hbar}
\Hbar_t = \frac{1}{2} H + \frac{1}{2} H b^{-1} \phi_t^{-1}.
\end{equation}
If $H$ is normalized, so is $\Hbar_t$.

If $b$ can be connected to the identity through a path 
of symplectomorphisms, then
the loop \eqref{shorter loop} is a deformation of the loop \eqref{phikt}.
Since $\Hbar_t(x) \leq \max H$ and $\Hbar_t(x) \geq \min H$
for all $x \in M$ and all $t$, 
\begin{equation} \labell{inequ}
\length( \{ \psi_t^{\Hbar} \} ) \leq \length ( \{ \psi_t^H \} ),
\end{equation}
and this is a strict inequality if $\max \Hbar_t < \max H$ for some $t$.

Denote by $\Nmax^H$ the subset of $M$ where $H$ takes its maximal value.
If $\Bmax^H \cap b(\Bmax^H) = 0$ then, for $t=0$,
the maximum of $\Hbar_0 = \half H + \half H b^{-1}$
is strictly smaller than that of $H$.

\smallskip

We have shown that,
if there exists a symplectic isotopy
that displaces the set where the momentum map is maximal,
then the non-effective circle action can be deformed to a loop 
of shorter Hofer length.
A similar argument gives a quantitative result:
if there exists a symplectic isotopy
that displaces the set where the momentum map is $\eps$-close to its maximum,
then the non-effective circle action can be deformed to a loop
whose Hofer length is shorter by $\eps$
from that of the original circle action.

\section{Shortening a combined action}
\labell{sec:polterov2}

In this section we give variations of Polterovich's trick
that apply to effective circle actions.
We begin with a lemma that we will use
to prove Theorem~\ref{thm:isolated}.
The lemma has two parts -- a qualitative part
and a quantitative part.

We say that a symplectomorphism $b$ \emph{disjoins}
a set $A$ from a set $B$ if $b(A) \cap B = \emptyset$.

{

\begin{Lemma} \labell{lem polterov}
Let $(M,\omega)$ be a compact symplectic manifold.
Let $H \colon M \to \R$ be the momentum map 
for a circle action $\{ \psi_t^H \}_{0 \leq t \leq 1}$.  
Let 
$$ i \colon (U,\omega_0) \to (M,\omega) $$
be an open symplectic embedding.
Suppose that
\begin{equation} \labell{HKF}
i^* H = K + F  \qquad \text{on \ $U$},
\end{equation}
where $K$ and $F$ generate commuting circle actions
$\{ \psi_t^{K} \}_{0 \leq t \leq 1}$ and 
$\{ \psi_t^{F} \}_{0 \leq t \leq 1}$
on $U$.   

Let $\hatN \colon U \to \R$ be a function 
that Poisson commutes with $K$ and $F$
and satisfies
\begin{equation}\labell{inf along N}
\inf\limits_{\hatN\inv(\alpha)} (K+F) = 
\inf\limits_{\hatN\inv(\alpha)} K + \inf\limits_{\hatN\inv(\alpha)} F 
\qquad \text{ for every $\alpha \in \R$. }
\end{equation}
Also suppose that 
\begin{equation} \labell{avoid min}
\text{ $i(U)$ does not meet the set where $H$ is minimal. }
\end{equation}

\begin{enumerate}
\item
Let 
\begin{align*}
\Bmax^H &= \{ m \in M \ | \ H \text{ is maximal at $m$ } \} , \\
\Bmax^K &= \{ u \in U \ | \ K \text{ is maximal at $u$ } \} , \\
\text{and} \qquad \Bmax^F &= \{ u \in U \ | \ F \text{ is maximal at $u$ } \} .
\end{align*}
Suppose that 
\begin{equation} \labell{included1}
\Bmax^H \subset i(U) \quad \text{ and } \quad 
\Bmax^K \cap \Bmax^F \neq \emptyset .
\end{equation}
Suppose that there exists a symplectomorphism $b \colon U \to U$
that disjoins $\Bmax^F$ from $\Bmax^K$
and such that $b$ can be connected to the identity through a path of 
symplectomorphisms that are compactly supported in $U$
and that preserve the function $\hatN$.

Then the circle action $\{ \Psi^H_t \}$ 
can be deformed, through loops in $\Ham(M,\omega)$,
to a loop $\{ \Psi^\Hbar_t \}$ whose Hofer length is smaller.

\item
Let $\eps > 0$.
Let
\begin{align*}
\Neps^H &= \{ m \in M \ | \ H(m) > \max H - \eps\}, \\
\Neps^K &= \{ u \in U \ | \  K(u) > \max K - \eps \}, \\
\text{and} \qquad 
\Neps^F &= \{ u \in U \ | \ F(u) > \max F - \eps \}.
\end{align*}
Suppose that 
\begin{equation} \labell{included2}
\Neps^H \subset i(U) \quad \text{ and } \quad 
\Bmax^K \cap \Bmax^F \neq \emptyset .
\end{equation}
Suppose that there exists a symplectomorphism $b \colon U \to U$
that disjoins $\Neps^F$ from $\Neps^K$
and such that $b$ can be connected to the identity through a path of 
symplectomorphisms that are compactly supported in $U$
and that preserve the function $\hatN$.

Then the circle action $\{ \Psi^H_t \}$ 
can be deformed, through loops in $\Ham(M,\omega)$,
to a loop $\{ \Psi^\Hbar_t \}$ whose Hofer length
is smaller than that of $\{ \Psi^H_t \}$ by at least $\eps$.
\end{enumerate}
In each of these two cases, the deformation can be chosen
such that the normalized deformed Hamiltonians 
remain between $\inf H(i(U))$ and $\max H(M)$
on the set $i(U)$ and coincide with $H$ outside this set,
%
%
so the positive Hofer length never exceeds the initial one
and the negative Hofer length remains the same throughout the deformation.
\end{Lemma}

\begin{proof}
Because $K$ and $F$ generate commuting circle actions,
$K \circ \psi_t^F = K$ and $F \circ \psi_t^K = F$ for all $t$.

Because $\Bmax^H \subset i(U)$,
\ $i^* H = K + F$, and $\Bmax^K \cap \Bmax^F \neq \emptyset$,
\begin{equation} \labell{maxHKF}
\Bmax^H = i \left( \Bmax^K \cap \Bmax^F \right) .
\end{equation}

Let $b \colon U \to U$ be a compactly supported symplectomorphism of $U$.

The loop
$$ \psi_t^{\Hbar} = 
\begin{cases}
i \circ \left( \psi_t^K \circ b \psi_t^F b\inv \right) \circ i\inv 
     & \text{ on } i(U) \\
\psi_t^H & \text{ outside } i(U)
\end{cases} $$
is generated by the Hamiltonian  $\Hbar_t \colon M \to \RR$
that is given by 
\begin{equation} \labell{Hbar with K and F}
\Hbar_t = 
\begin{cases}
    \left( K + F b\inv (\psi_t^K)\inv \right) \circ i\inv 
      & \text{ on } i(U) \\
    H & \text{ outside } i(U).
\end{cases} 
\end{equation}

If $H$ is normalized, so is $\Hbar_t$.

If $b$ can be connected to the identity through 
a compactly supported symplectic isotopy,
then the loop $ \{ \psi^\Hbar_t \}$
is a deformation of the loop $ \{ \psi^H_t \}$
in $\Ham(M,\omega)$.

By~\eqref{maxHKF} and~\eqref{Hbar with K and F},
we have
$\max \Hbar_t \leq \max H$ for all $t$. So the positive Hofer length
of $\{ \psi^\Hbar_t \}$ does not exceed that of the 
circle action $\{ \psi^H_t \} $.

If $b$ preserves $\hatN$, then, 
by~\eqref{inf along N}, \eqref{avoid min}, \eqref{Hbar with K and F},
and because $\psi_t^K$ preserves $\hatN$,
we have $\min \Hbar_t = \min H$ for all $t$. 
So the negative Hofer length of $\{ \psi^\Hbar_t \}$
is equal to that of the circle action $\{ \psi^H_t \}$.

Suppose that $b$ disjoins $\Bmax^F$ from $\Bmax^K$.
Then, for $t=0$, we have
$\max \Hbar_0 = \max (K + Fb\inv) < \max K + \max F = \max H$.
This strict inequality is because
$ \Bmax^K \cap \Bmax^{Fb\inv} $,
being equal to $\Bmax^K \cap b(\Bmax^F)$, is empty.
(In fact, $\max \Hbar_t < \max H$ for \emph{all} $t$;
this follows from 
$K + F b\inv (\psi^K_t)\inv = (K+Fb\inv) (\psi^K_t)\inv$
and $\Bmax^K \cap b(\Bmax^F) = \emptyset$.)
Because the inequality $\max \Hbar_t \leq \max H$ 
holds for all $t$ and is strict for $t=0$,
the positive Hofer length of $\{ \psi^\Hbar_t \}$ 
is strictly smaller than that of the circle action $\{ \psi^H_t \}$.

This completes the proof of the first, qualitative, 
part of the lemma.

To prove the second, quantitative, part of the lemma,
suppose now that $b$ disjoins $\Neps^F$ from $\Neps^K$,
let $m$ be a point in $M$,
and we will show that $\Hbar_t(m) \leq \max H - \eps$ for all $t$.

Suppose that $m$ is not in $i(U)$, 
and, thus, also not in $\Neps^H$.
Then $\Hbar_t(m) = H(m)$ and $H(m) \leq \max H - \eps$, so
$$ \Hbar_t(m) \leq \max H - \eps.$$
Now suppose that $m$ is in $i(U)$, say, $m=i(u)$ for $u \in U$.  
By~\eqref{Hbar with K and F},
$$ \Hbar_t(m) = K(u) + F(u')$$ 
where $u' = b\inv (\psi^K_t)\inv (u)$.
If $\Hbar_t(m)$ is $\eps$-close to $\max H$,
then $u \in \Neps^K$ and $u' \in \Neps^F$, that is,
\begin{equation} \labell{m}
u \in \Neps^K \cap \psi^K_t(b(\Neps^F)) .
\end{equation}
Because $K = K \circ \psi^K_t$, 
The set $\Neps^K$ is preserved under $\psi^K_t$,
so we can rewrite~\eqref{m} as 
$u \in \psi^K_t(\Neps^K \cap b(\Neps^F))$.
But, by assumption, $\Neps^K \cap b(\Neps^F)$ is empty.
So $\Hbar_t(m)$ cannot be $\eps$-close to $\max H$,
and again we conclude that $\Hbar_t(m) \leq \max H - \eps$. 
So $\Hbar_t(m) \leq \max H - \eps$ for all $m \in M$
and for all $t$.
Since $\Hbar_t$ is normalized, we conclude that
the positive Hofer lengths satisfy
$\ell_+ ( \{ \psi_t^{\Hbar} \} ) \leq \ell_+ ( \{ \psi_t^H \} ) - \eps$, 
as required. 
\end{proof}

}

The second, quantitative, part of Lemma~\ref{lem polterov}
does not automatically extend to more than two summands.
The first, qualitative, part does extend to more than two summands.
We now give such an extension.  We will use it
to prove Theorems~\ref{thm:nonisolated} and~\ref{shorten w subbundle}.

\begin{Lemma} \labell{lem polterov k}
Let $(M,\omega)$ be a compact symplectic manifold.
Let $H \colon M \to \R$ be the momentum map for a circle action
$\{ \psi^H_t \}_{0 \leq t \leq 1}$.  Let
$$ i \colon (U,\omega_0) \to (M,\omega) $$ 
be an open symplectic embedding.
Suppose that
\begin{equation} \labell{iH}
i^*H = H_{(0)} + H_{(1)} + \ldots + H_{(k)} \qquad \text{ on $U$}
\end{equation}
where $H_{(0)}, H_{(1)}, \ldots, H_{(k)}$ generate commuting circle actions
on $U$.  

Let
$$ \Bmax^H = \{ m \in M \ | \ H \text{ is maximal at $m$} \} $$
and, for each $0 \leq j \leq k$,
$$ \Bmax^{(j)} = \{ u \in U \ | \ H_{(j)} \text{ is maximal at $u$} \}.$$
Suppose that
\begin{equation} \labell{included}
\Bmax^H \subset i(U) \quad \text{ and } \quad
\Bmax^{(0)} \cap \Bmax^{(1)} \cap \ldots \cap \Bmax^{(k)} \neq \emptyset .
\end{equation}

Suppose that 
\begin{equation} \labell{inf along N again}
\inf\limits_{U} \left( H_{(0)} + H_{(1)} + \ldots + H_{(k)} \right)  
   =   \sum_{j=0}^k \inf\limits_{U} H_{(j)}
\end{equation}
and that 
\begin{equation} \labell{avoid min again}
\text{ $i(U)$ does not meet the set where $H$ is minimal. } 
\end{equation}

Suppose that there exist symplectomorphisms $b_1,\ldots,b_k$ of $U$
such that 
\begin{equation} \labell{k disjoin}
\Bmax^{(0)} \cap b_1(\Bmax^{(1)}) \cap \ldots \cap b_k(\Bmax^{(k)}) 
= \emptyset 
\end{equation}
and such that each $b_j$ 
can be connected to the identity through a path of symplectomorphisms
that is compactly supported in $U$. 

Then the circle action $\{ \psi^H_t \}$ can be deformed,
through loops in $\Ham(M,\omega)$, to a loop $\{\psibar_t\}$
whose Hofer length is smaller.

The deformation can be chosen such that the normalized deformed Hamiltonians 
remain between $\inf H(i(U))$ and $\max H(M)$ on the set $i(U)$
and coincide with $H$ outside this set,
so the positive Hofer length never exceeds the initial one 
and the negative Hofer length remains the same.
\end{Lemma}

\begin{proof}
Denote by
$$ \{ \psi^{(0)}_t \}_{0 \leq t \leq 1} \ , \ \ldots \ , \  
\{ \psi^{(k)}_t \}_{0 \leq t \leq 1} $$
the circle actions on $U$ that are generated by $H_{(0)},\ldots,H_{(k)}$.
Let $b_1,\ldots,b_k$ be arbitrary compactly supported symplectomorphisms
of $U$.  For each $j \in \{ 0 , \ldots , k \}$, let
\begin{equation} \labell{psibar j}
\psibar^{(j)}_t: = \psi^{(0)}_t \circ \left(b_1 \psi^{(1)}_t b_1\inv \right)
\circ \cdots \circ \left( b_j \psi^{(j)}_t b_j\inv \right) .
\end{equation}
This is a loop in $\Ham(U,\omega_0)$.  
We have the following three facts:
\begin{enumerate}
\item
$\{ \psibar^{(j)}_t \}$ is generated by the Hamiltonian
\begin{equation} \labell{Hbar j t}
\Hbar^{(j)}_t := H_{(0)} + H_{(1)} b_1\inv (\psibar^{(0)}_t)\inv 
+ \ldots + H_{(j)} b_j\inv (\psibar^{(j-1)}_t)\inv .
\end{equation}
\item
Outside a compact subset of $U$,
$$ \psibar^{(j)}_t = 
\psi^{(0)}_t \circ \psi^{(1)}_t \circ \ldots \circ \psi^{(j)}_t
\quad \text{ and } \quad
\Hbar^{(j)}_t = H_{(0)} + H_{(1)} + \ldots + H_{(j)} .$$
\item
For all $t$,
\begin{equation} \labell{3}
\int_U \left( \Hbar^{(j)}_t - (H_{(0)} + H_{(1)} + \ldots + H_{(j)} ) 
  \right) \omega_0^n = 0 .
\end{equation}
\end{enumerate}
Facts (1) and (2) are shown by induction on $j$.
As for Fact (3), by~\eqref{Hbar j t},
the left hand side of~\eqref{3} is the sum over $m=1,\ldots,j$ of
\begin{equation} \labell{annie}
\int_U H_{(m)} 
  \circ \left( b_{m} \inv \circ (\psibar^{(m-1)}_t)\inv \right) \omega_0^n 
- \int_U H_{(m)} \omega_0^n .
\end{equation}
The map $b_m\inv \circ ( \psibar^{(m-1)}_t )\inv$ is a symplectomorphism, 
so it preserves the volume form $\omega_0^n$, 
and so~\eqref{annie} vanishes.


By the above items (1) and (2), we can define a loop in $\Ham(M,\omega)$ by
$$ \psibar_t := \begin{cases}
i \circ \psibar^{(k)}_t \circ i\inv 
			     & \text{ on } \ i(U) \\
\psi^H_t & \text{ on } \ M \ssminus i(U) ,
\end{cases} $$
and it is generated by the Hamiltonian
%
%
%
\begin{equation} \labell{Hbar t}
\Hbar_t := \begin{cases}
\Hbar^{(k)}_t \circ i\inv & \text{ on } \ i(U) , \\
H & \text{ on } \ M \ssminus i(U).
\end{cases}
\end{equation}

By the above item (3), if $H$ is normalized, so is $\Hbar_t$.  

From~\eqref{iH}, 
\eqref{inf along N again}, \eqref{avoid min again}, 
\eqref{Hbar j t}, and~\eqref{Hbar t}, 
we deduce that $\min\limits_M \Hbar_t = \min\limits_M H$.
So the negative Hofer length of $\{ \psibar_t \}$
is equal to that of the circle action~$\{ \psi^H_t \}$.

From~\eqref{iH}, \eqref{included}, 
                \eqref{Hbar j t}, 
and~\eqref{Hbar t}, 
we see that $\max\limits_M \Hbar_t \leq \max\limits_M H$ for all $t$,
with a strict inequality when $t=0$ 
if $b_1,\ldots,b_k$ satisfy~\eqref{k disjoin}.
Integrating, we get the inequality 
$\int_0^1 \max\limits_M \Hbar_t dt \leq \int_0^1 \max\limits_M H dt$,
and if $b_1,\ldots,b_k$ satisfy~\eqref{k disjoin}
then this inequality is strict
(because in this case the pointwise inequality between the integrands
is strict in a neighbourhood of $t=0$).
That is, the positive Hofer length of $\{ \psibar_t \}$
does not exceed that of $\{ \psi^H_t \}$,
and, if $b_1,\ldots,b_k$ satisfy~\eqref{k disjoin},
then the positive Hofer length of $\{ \psibar_t \}$ is strictly smaller
than that of $\{ \psi^H_t \}$.

The lemma is obtained by applying this argument
as $b_1,\ldots,b_k$ vary over smooth paths of symplectomorphisms 
that start with the identity,
are compactly supported in $U$,
and end with $b_1,\ldots,b_k$ that satisfy~\eqref{k disjoin}.


\end{proof}

\section{Disjoining discs in $\CC$}
\labell{sec:disjoin disc}

For every neighbourhood
of the origin in $\C$, there exists a Hamiltonian
isotopy, supported in this neighbourhood,
that disjoins the origin from itself.
Such an isotopy can be generated
by a function $F \colon \C \to \R$
whose Hamiltonian vector field
is zero outside the given neighbourhood,
is equal to $-\dd{x}$ on a neighbourhood of the origin,
and is equal to a non-negative multiple of $-\dd{x}$
everywhere on the $x$-axis.
For example, if the given neighbourhood contains the set
$\{z \ | \ \pi|z|^2 \leq 3 \eps \}$,
we may take the function $F(z) = y \, \rho(\pi|z|^2/\eps - 1)$
where $y$ is the imaginary part of $z$
and where $\rho(s)$ is equal to one for $s \leq 0$,
to zero for $s \geq 1$, and to
$(\int_s^1 e^{-\frac{1}{\tau(1-\tau)}} d\tau)
/(\int_0^1 e^{-\frac{1}{\tau(1-\tau)}} d\tau) $
for $0 \leq s \leq 1$.

For our quantitative results we will need to
keep track of the sizes of neighbourhood
of the origin that get disjoined.
We do so in the following lemma.

\begin{Lemma} \labell{lem:btA}
There exists a smooth family of functions
$$ F^{A,\eps}_t \colon \C \to \R \quad , \quad 0 \leq t \leq 1 ,$$ 
parametrized by $A \geq 0$ and $0 < \eps < 1$, 
such that $F^{A,\eps}_t(z) = 0$ for $z$ outside the disc
$ \{ z \ | \ \pi |z|^2 < 1 + A + \eps \} $
and such that the Hamiltonian flow $b^{A,\eps}_t \colon \C \to \C$ 
of $F^{A,\eps}_t$, at time $t=1$, carries the disc
$ \{ z \ | \ \pi|z|^2 \leq 1 \} $
into the annulus
$ \{ z \ | \ A < \pi|z|^2 < 1 + A + \eps \} $.
\end{Lemma}

\begin{Remark}
\emph{Smooth family} means that the function
$(A,\eps,t,z) \mapsto F_t^{A,\eps}(z)$
from $[0,\infty) \times (0,1) \times [0,1] \times \C$
to $\R$ is smooth in the usual sense.
It follows that the function
$(A,\eps,t,z) \mapsto b_t^{A,\eps}(z)$ is smooth.
\end{Remark}

Lemma~\ref{lem:btA} is used in the next section,
in the proof of Lemma~\ref{CtimesW}.
The rest of this section is devoted to the proof
of Lemma~\ref{lem:btA}.  The reader may skip the proof
and proceed to the next section.

\begin{Lemma} \labell{lem:bt1}
There exists a smooth family of (time independent) functions 
$H^\eps \colon \C \to \R$, parametrized by $0 < \eps < 1$,
whose Hamiltonian flows $b^\eps_t$
have the following properties.
\begin{itemize}
\item
For every $t$, the diffeomorphism $b^\eps_t$ 
is the identity map outside the disc 
\begin{equation} \labell{disc1pluseps}
\{ z \ | \ \pi|z|^2 < 1+\eps \}.
\end{equation}
\item
For $t>1$, the diffeomorphism $b^\eps_t$ carries 
the disc 
\begin{equation} \labell{disc1}
\{ z \ | \ \pi|z|^2 \leq 1 \}
\end{equation}
into the slit disc 
\begin{equation}  \labell{slit disc}
\{ z \ | \ \pi|z|^2 < 1 + \eps \} 
\ssminus \text{the non-negative x-axis} .
\end{equation}
\end{itemize}
\end{Lemma}

\begin{proof}
We choose $H^\eps$ such that, for each $\eps$, 
the Hamiltonian vector field of $H^\eps$ is zero 
outside the disc~\eqref{disc1pluseps},
is equal to $-\frac{1}{\sqrt{\pi}}\dd{x}$ on the segment
$[0,\frac{1}{\sqrt{\pi}}]$ of the x-axis,
and is equal to a non-negative multiple of $-\dd{x}$ 
on the rest of the x-axis.

For example, we may take 
$H^\eps(z) = \frac{1}{\sqrt{\pi}} y \rho( (\pi|z|^2-1)/\eps)$
where $y$ is the imaginary part of $z$ and where
$\rho(s)$ is equal to $1$ for $s \leq 0$,
to $0$ for $s \geq 1$, and is non-negative.

For each $t$,
the Hamiltonian flow generated by the function $H^\eps$,
at times $t \geq 0$, carries the segment 
$[-\sqrt{\frac{1+\eps}{\pi}},\frac{1}{\sqrt{\pi}}]$ 
to a segment of the form $[-\sqrt{\frac{1+\eps}{\pi}},x_t]$;
if $t>1$, then $x_t<0$.
Also, this flow carries the segment 
$[\frac{1}{\sqrt{\pi}},\sqrt{\frac{1+\eps}{\pi}}]$ 
to the segment $[x_t,\sqrt{\frac{1+\eps}{\pi}}]$,
and it is the identity map outside the disc~\eqref{disc1pluseps}.
Thus, it satisfies the requirements of the lemma.
\end{proof}

\begin{Lemma} \labell{lem:bt2}
There exists a smooth family of functions
$H^{\eps,\delta}_t \colon \C \to \R$,
parametrized by $0<\eps<1$ and $0 < \delta < \eps$,
whose Hamiltonian flows $b^{\eps,\delta}_t \colon \C \to \C$
have the following properties.
\begin{itemize}
\item
The diffeomorphism $b^{\eps,\delta}_t$
is the identity map outside the disc
$\{ z \ | \ \pi|z|^2 < t+1+\eps \}$.
\item
The diffeomorphism $b^{\eps,\delta}_t$ carries the set
\begin{equation} \labell{minuswedge}
\{ re^{i\theta} \ | \ \delta \leq \pi r^2 \leq 1 + \eps - \delta
\ , \ \delta \leq \theta \leq 2\pi - \delta \} 
\end{equation}
to the set
\begin{equation} \labell{minuswedge plus t}
\{ re^{i\theta} \ | \ t+\delta \leq \pi r^2 \leq t + 1 + \eps - \delta
\ , \ \delta \leq \theta \leq 2\pi - \delta \} .
\end{equation}
\end{itemize}
\end{Lemma}

\begin{proof}
We choose $H^{\eps,\delta}_t$ whose Hamiltonian vector field
vanishes outside the annulus
$  \{ z \ | \ t < \pi |z|^2 < 1 + t + \eps \}$
and, on the set~\eqref{minuswedge plus t},
it is equal to $\frac{1}{2\pi r} \dd{r}$.
(Then, on the set~\eqref{minuswedge plus t},
the derivative of the function $\pi |z|^2$
along this vector field is equal to one.)
The Hamiltonian flow generated by the 
function~$H^{\eps,\delta}_t$, at time $t \geq 0$,
will carry the set~\eqref{minuswedge}
to the set~\eqref{minuswedge plus t}, as required.

To be explicit, we may take 
$H^{\eps,\delta}_t(re^{i\theta}) = 
- \frac{\theta}{2\pi}
\rho_1^{\delta}(\theta) \rho_2^{\eps,\delta}(\pi r^2 - t) $
where $0 \leq \theta < 2\pi$
and $t \leq \pi r^2 \leq 1+t+\eps$,
where $\rho_1^{\delta} \colon [0,2\pi] \to \R$
is a smooth family of functions that vanish to all orders
at $\theta=0$ and $\theta=2\pi$
and such that $\rho_1^\delta$ is equal to one on $[\delta,2\pi-\delta]$,
and where $\rho_2^{\eps,\delta} \colon \R \to \R$ is a smooth family 
of functions
that vanish outside $[0,1+\eps]$ and such that $\rho_2^{\eps,\delta}$
is equal to one on $[\delta,1+\eps-\delta]$.
\end{proof}

\begin{proof}[Proof of Lemma~\ref{lem:btA}]

Every closed subset of $\C$ which is contained in the slit disc
\begin{equation} \labell{slit disc again}
\{ z \ | \ \pi|z|^2 < 1 + \eps \}  
\ssminus \text{the non-negative x-axis} 
\end{equation}
is also contained in a set of the form 
$$ \{ re^{i\theta} \ | \ \delta \leq \pi r^2 \leq 1 + \eps - \delta
\ , \ \delta \leq \theta \leq 2\pi - \delta \} $$
for some $\delta > 0$.
Moreover, for every closed subset of $(0,1) \times \C$
which is contained in the product of the open segment $(0,1)$
with the slit disc~\eqref{slit disc again},
there exists a positive smooth function $\eps \mapsto \delta_\eps$
such that the closed subset is also contained in the subset
\begin{equation} \labell{delta eps}
\{ (\eps, re^{i\theta}) \ | \ 
0 < \eps < 1 \ , \ 
\delta_\eps \leq \pi r^2 \leq 1 + \eps - \delta_\eps
\ , \ \delta_\eps \leq \theta \leq 2\pi - \delta_\eps \} 
\end{equation}
of $(0,1) \times \C$.

Let $b^\eps_t$ be the Hamiltonian flows of Lemma~\ref{lem:bt1}.
Because $(\eps,z) \mapsto (\eps,b^\eps_t(z))$
is a diffeomorphism of $(0,1) \times \C$,
it carries the set 
$\{(\eps,z) \ | \ \pi|z|^2 \leq 1 \}$
to a closed subset of $(0,1) \times \C$.
For $t>1$, by the second item of Lemma~\ref{lem:bt1},
this closed subset is contained in the product
of $(0,1)$ with the slit disc~\eqref{slit disc again}.
Fix $T > 1$, and fix a smooth function $\eps \mapsto \delta_\eps$
such that the set
$\{ (\eps,b^\eps_T(z)) \ | \ \pi|z|^2 \leq 1 \}$
is contained in the set~\eqref{delta eps}.

Let $b^{\eps,\delta}_t$ be the Hamiltonian flows
of Lemma~\ref{lem:bt2}.

Let $\rho \colon [0,1] \to [0,1]$ be a smooth function
that takes $0$ to $0$ and $1$ to $1$
and whose derivatives of all orders vanish at the endpoints 
$0$ and $1$.  For example, we may take
$\rho(s) = \int_0^s e^{-\frac{1}{s(1-s)}}ds / 
   \int_0^1 e^{-\frac{1}{s(1-s)}}ds $.
The Hamiltonian flow
$$ b^{A,\eps}_\tau = \begin{cases}
b^\eps_{T \rho(2\tau)} & 0 \leq \tau \leq 1/2 \\
b^{\eps,\delta_\eps}_{A\rho(2\tau-1)} & 1/2 \leq \tau \leq 1 
\end{cases}$$
satisfies the requirements of Lemma~\ref{lem:btA}.
\end{proof}

\section{Disjoining a family of discs}
\labell{sec:disjoin family}

The purpose of this section is to prove the following
quantitative result, which is later used in the proof
of Theorem~\ref{thm:isolated}.
Recall that a function from a topological space $X$
to a topological space $Y$ is \emph{proper}
if the preimage of every compact subset of $Y$
is a compact subset of $X$.

\begin{Lemma} \labell{CtimesW}
Let $W$ be a symplectic manifold and $N \colon W \to [0,r)$
a proper function.  Let $A_1,A_2 \colon [0,r) \to \R$
be smooth functions such that the set
$$ C := \{ x \in [0,r) \ | \ A_1(x) \geq 0 \text{ and } A_2(x) \geq 0 \} $$
is compact.  Let $U$ be an open subset of $\C \times W$ that contains
the set
\begin{equation} \labell{subset of CtimesW}
\left\{ (z,w) \ | \ N(w) \in C \text{ and } 
	     \pi |z|^2 \leq A_1(N(w)) + A_2(N(w)) \right\}.
\end{equation}
Then there exists a symplectic isotopy, compactly supported in $U$,
that disjoins the set
$$ \left\{ (z,w) \ | \ \pi|z|^2 \leq A_1(N(w)) \right\} $$
from the set
$$ \left\{ (z,w) \ | \ \pi|z|^2 \leq A_2(N(w)) \right\} $$
and that preserves the function $(z,w) \mapsto N(w)$.
\end{Lemma}

Lemma~\ref{CtimesW} is, more or less, a parametrized version
of the fact that a disc of area $A_1$ can be disjoined
from a disc of area $A_2$ inside any disc of area greater than $A_1+A_2$.
The symplectic isotopy moves the $w$ coordinate but it does not change 
the value of $N(w)$.  For a fixed value of $N(w)$, the effect of the
isotopy on the $z$ coordinate is independent of the $w$ coordinate.

\begin{proof}
The map $(z,w) \mapsto \left( \pi|z|^2 , N(w) \right)$
from $\C \times W$ to $[0,\infty) \times [0,r)$ is proper, hence closed.  
So the image of the complement of $U$ under this map is a closed
subset of $[0,\infty) \times [0,r)$.
This image is disjoint from the compact set 
$$ \{ (s,x) \in [0,\infty) \times [0,r) \ | \ x \in C \text{ and }
s \leq A_1(x) + A_2(x) \}$$ 
so it has a positive distance from this set.
Let $\epsbar$ be a positive number that is smaller than half 
of this distance
and is smaller from the difference $r - \max C$.
Then $V := ( C + (-\epsbar,\epsbar) ) \cap [0,r)$ is an open 
neighbourhood of $C$ in $[0,r)$ whose closure 
in $[0,r)$ is compact, and $U$ contains the set
\begin{equation} \labell{ora}
\left \{ (z,w) \ | \ N(w) \in \text{closure}(V) \text{ and }
\pi|z|^2 \leq A_1(N(w)) + A_2(N(w)) + \epsbar  \right\} .
\end{equation}


Let $\Fbar^{A_1,A_2}_t \colon \C \to \R$ be a smooth family of functions,
defined for $t \in [0,1]$ and $A_1,A_2 \in \R$,
such that, when $A_1 \geq 0$ and $A_2 \geq 0$, 
the time-dependent function $z \mapsto \Fbar_t^{A_1,A_2}(z)$
vanishes outside the set given by $\pi |z|^2 \leq A_1 + A_2 + \epsbar$
and its Hamiltonian flow $\bbar^{A_1,A_2}_t$, at time $t=1$, 
carries the set $\{ z \ | \ \pi |z|^2 \leq A_1 \}$
into the set $\{ z \ | \ A_2 < \pi |z|^2 \leq A_1 + A_2 + \epsbar \}$.
For example, we may set
$$ \Fbar^{A_1,A_2}_t(z) := 
c F^{A,\eps}_t \, (z/c)$$
where $c = \sqrt{A_1+\epsbar/2}$, 
$A=\displaystyle{\frac{A_2}{A_1+\epsbar/2}}$,
and $\eps=\displaystyle{\frac{\epsbar/2}{A_1+\epsbar/2}}$,
and where $F^{A,\eps}_t$ is as in Lemma~\ref{lem:btA}.
(The Hamiltonian flow $\bbar^{A_1,A_2}_t$ of $\Fbar^{A_1,A_2}_t$
relates to the Hamiltonian flow $b^{A,\eps}_t$ of $F^{A,\eps}_t(z)$
by $\bbar^{A_1,A_2}_t(z) = c b^{A,\eps}_t(z/c)$.)

Let $\rho \colon [0,r) \to \R$ be a smooth function 
that vanishes outside $V$ and such that $\rho|_C \equiv 1$.
Let 
$$ H_t(z,w) = \rho(N(w)) \Fbar^{A_1(N(w)),A_2(N(w))}_t (z) .$$
This function vanishes outside the compact set~\eqref{ora}.
Its Hamiltonian flow preserves the function $N(w)$.
We will show, for every~$N_0$, that this flow, at time $t=1$,
disjoins 
\begin{equation} \labell{setA1}
\{ (z,w) \ | \ \pi|z|^2 \leq A_1(N_0) \text{ and } N(w) = N_0 \}
\end{equation}
from
\begin{equation} \labell{setA2}
\{ (z,w) \ | \ \pi|z|^2 \leq A_2(N_0) \text{ and } N(w) = N_0 \}.
\end{equation}

When $N_0 \in C$, the restriction of the Hamiltonian flow of $H_t$
to the level set $\{ N(w) = N_0 \}$ is
$$ b_t(z,w) = \left( \bbar^{A_1,A_2}_t , b^N_{T(z,w,t)}(w) \right) $$
where $A_1=A_1(N_0)$, where $A_2=A_2(N_0)$,
where $\bbar^{A_1,A_2}_t$ is the Hamiltonian flow
of $\Fbar^{A_1,A_2}_t$ on $\C$,
where $b^N_T$ is the Hamiltonian flow of $N$ on $W$
(with time parameter $T$),
and where $T(z,w,t)$ is some real valued function of $z,w,t$.
At $t=1$, this flow carries the set 
$$ \left\{ (z,w) \ | \ \pi|z|^2 \leq A_1(N_0) 
\text{ and } N(w) = N_0 \right\} $$
into the set
$$ \left\{ (z,w) \ | \ 
    A_2(N_0) < \pi|z|^2 \leq A_1(N_0) + A_2(N_0) + \epsbar
    \text{ and } N(w) = N_0 \right\},$$
so it disjoins~\eqref{setA1} from~\eqref{setA2}.

When $N_0 \not\in C$, one of the sets~\eqref{setA1}
and~\eqref{setA2} is empty.
So the flow trivially disjoins 
the set~\eqref{setA1} from the set~\eqref{setA2}. 
\end{proof}

\section{Shortening on a manifold with an isolated maximum}
\labell{sec:isolated}

In this section we prove Theorem~\ref{thm:isolated}
and Remarks~\ref{rk:isolated} and~\ref{rk:isolated'}.
We recall the statement:

\begin{quotation}
Let $(M,\omega)$ be a compact connected symplectic manifold
with a Hamiltonian circle action.
Suppose that the momentum map attains its maximum 
at an isolated fixed point
and that at least one of the isotropy weights at that point
has absolute value greater than one.
Let~$d$ be a positive number that is smaller than the gap
between the two largest critical values of the momentum map.
Then the circle action can be deformed to a loop in $\Ham(M,\omega)$ 
whose Hofer length is smaller than that of the original circle action
by $2d/9$.

We can choose the deformation
such that the positive Hofer length never exceeds the initial one
and the negative Hofer length remains constant.
Specifically, let $H$ denote the momentum map, let $\pmax$ denote
the point where $H$ attains its maximum, let $\alpha$ denote the
gap between the two largest critical values of $H$,
and consider the subset of $M$ given by
\begin{equation} \labell{alpha nbhd}
\{ m \ | \ H(m) > H(\pmax) - \alpha \}.
\end{equation}
We can choose the deformation such that 
the normalized deformed Hamiltonians remain greater than $H(\pmax)-\alpha$
on the set~\eqref{alpha nbhd} and coincide with $H$ outside this set.

We actually obtain a better estimate than $\frac{2}{9}d$.
If one of the isotropy weights is even then we can shorten by $\frac{1}{4}d$.
Otherwise, let $-k_1$ be the weight whose absolute value is largest;
then we can shorten by $\frac{1}{4} ( 1 - \frac{1}{k_1^2}) d$.
When $k_1=3$, this becomes $\frac{2}{9}d$.
\end{quotation}

\begin{proof}[Proof of Theorem~\ref{thm:isolated}
and Remarks~\ref{rk:isolated} and~\ref{rk:isolated'}]
Let $-k_1,\ldots,-k_n$ be the isotropy weights at $\pmax$.
The equivariant Darboux theorem implies that the integers $k_j$ 
are positive.  We assume that at least one of them 
is greater than one.
Without loss of generality,
\begin{equation} \labell{kab}
k_1 = a + b 
\end{equation}
where $a$ and $b$ are positive integers.

The point $\pmax$ is the only critical point of $H$ in the 
set~\eqref{alpha nbhd},
and the momentum map $H$ is proper as a map from this set
to the ray $(H(\pmax)-\alpha,\infty)$.
There exists an equivariant symplectomorphism
from the subset of $\C^n$ given by  
$$\calN^H_\alpha 
= \{ z \ | \ 
 \pi(k_1 |z_1|^2 + \ldots + k_n |z_n|^2) < \alpha \},$$
with the circle action generated by the momentum map
$H(z_1,\ldots,z_n) = -k_1 \pi |z_1|^2 - \ldots - k_n |z_n|^2$
to the subset~\eqref{alpha nbhd} of $M$;
denote it
$$ i \colon \calN^H_\alpha \to M .$$
This follows from \cite[Proposition~2.8]{KT}
(which should have been stated with $|z_j|^2$ instead of the typo $|z_j|$).

Let
$$ N(z_2,\ldots,z_n) := \pi (k_2|z_2|^2 + \ldots + k_n|z_n|^2).$$

We will apply Lemma~\ref{lem polterov} with $U = \calN^H_\alpha$,
$$ K(z_1,\ldots,z_n) := H(\pmax) - \pi b |z_1|^2 ,$$
and
$$ F(z_1,\ldots,z_n) := - \pi a |z_1|^2 - N(z_2,\ldots,z_n) .$$ 

Let $0 < d < \alpha$ and let $s = \frac{2}{9} d$.
By Lemma~\ref{lem polterov},
to shorten by amount $s$, it is enough to disjoin the set 
$\calN_s^K = \{ z \ | \ \pi b |z_1|^2 < s \}$
from the set $\calN_s^F = 
\{ z \ | \ \pi a |z_1|^2 + N(z_2,\ldots,z_n) < s \}$
by a compactly supported symplectic isotopy of ${\calN}^H_\alpha$
that preserves the function 
$$\hatN(z_1,z_2,\ldots,z_n) := N(z_2,\ldots,z_n).$$

Write the sets $\calN_s^K$ and $\calN_s^F$ as 
$$\calN_s^K = \{ z \ | \ \pi |z_1|^2 < \frac{s}{b} \}
\quad \text{and} \quad
\calN_s^F = \{ z \ | \ \pi |z_1|^2 < \frac{s - N}{a} \}.$$
By Lemma~\ref{CtimesW}, with $A_1(N) = s/b$ and $A_2(N) = (s-N)/a$,
for any open neighbourhood $U$ of the set 
\begin{equation} \labell{ambient}
\{ z \ | \ \pi |z_1|^2 \leq \frac{s}{b} + \frac{s - N}{a} \}
\end{equation}
there exists a symplectic isotopy, compactly supported in $U$,
that disjoins $\calN_s^F$ from $\calN_s^K$
and that preserves the function $\hatN$.
So it is enough to show that the set~\eqref{ambient}
is contained in the set $\calN_\alpha^H$.

Because $a+b=k_1$, we can rewrite the set~\eqref{ambient} as
\begin{equation} \labell{set1}
\{ z \ | \ \frac{ab}{sk_1} \pi |z_1|^2 + \frac{b}{sk_1} N \leq 1 \}.
\end{equation}
Also, we rewrite the set $\calN_\alpha^H$ as
\begin{equation} \labell{set2}
\{ z \ | \ \frac{k_1}{\alpha} \pi |z_1|^2 + \frac{1}{\alpha} N < 1 \}.
\end{equation}
It remains to show that the set~\eqref{set1} is contained 
in the set~\eqref{set2}.  This is equivalent to requiring that
$$ \frac{ab}{sk_1} \geq \frac{k_1}{\alpha}
\quad \text{and} \quad
\frac{b}{sk_1} \geq \frac{1}{\alpha}. $$
Because $a \leq k_1$, the first of these inequalities implies the second.

If $k_1$ is even, write $k_1 = a+b$ where $b=a$.
Then $ab/k_1^2 = 1/4$, which is greater than $2/9$.
If $k_1$ is odd, write $k_1 = a+b$ where $b=a+1$.
Then $ab/k_1^2 = 1/4 ( 1 - 1/k_1^2 )$
which, because $k_1 \geq 3$, is greater than or equal to 
$1/4 ( 1 - 1/9 ) = 2/9$.
In either case, if $s = 2/9 d$ and $d < \alpha$, then 
$s < (ab/k_1^2)\alpha$, which is equivalent to the required inequality.
\end{proof}



\section{Symplectic vector bundles}
\labell{sec:nbhd of B}

In this section we apply the symplectic tubular neighbourhood theorem
to describe a neighbourhood of a fixed component
in a manifold with Hamiltonian circle action in terms of a symplectic
vector bundle,
and we give some results about symplectic vector bundles
that we will use in the next section.


\begin{Lemma} \labell{omegaTheta}
Let $$\pi \colon E \to B$$ be a symplectic vector bundle,
i.e., a vector bundle with a fibrewise symplectic form.
There exists on the total space $E$ a closed two-form $\omega_\Theta$ 
with the following properties.
\begin{enumerate}
\item
The pullback of $\omega_\Theta$ to the fibres of $E$
coincides with the fibrewise symplectic forms.
\item
The pullback of $\omega_\Theta$ by the zero section of $E$
is zero.
\item
At the points of the zero section, the fibres of $E$
are $\omega_\Theta$-orthogonal to the zero section.
\end{enumerate}
Moreover, if a compact Lie group $K$ acts on $E$ by bundle 
automorphisms,
then $\omega_\Theta$ can be chosen to be $K$ invariant.
\end{Lemma}

\begin{proof}
Let $2n = \text{rank} E$, and let
$$ G = \Sp(\R^{2n}) $$  
be the group of symplectic linear transformations of $\R^{2n}$.
Let $P \to B$ be the principal $G$ bundle 
(with a right $G$ action) whose fibre $E_b$ over a point $b$ of $B$ 
consists of the set of linear symplectic isomorphisms 
from $\R^{2n}$ to $E_b$.

We identify $E$ with the associated bundle 
$$ P \times_G \R^{2n} , $$
in which $[pa,z] = [p,az]$ for all $p \in P$, $z \in \R^{2n}$,
and $a \in G$.

Let $\Theta \in \Omega^1(P,\g)$ be a connection one-form;
use the same symbol, $\Theta$, to denote its pullback
to $P \times \R^{2n}$.
Let $\omega_0$ be the standard symplectic form on $\R^{2n}$.
Let $\Phi_{\R^{2n}} \colon \R^{2n} \to \g^* $
be the quadratic momentum map for the $G$-action on $\R^{2n}$.
Let $\left< \cdot , \cdot \right>$ denote the pairing
between $\g^*$ and $\g$.  Then
$\left< \Phi_{\R^{2n}} , \Theta \right>$
is a $G$-invariant real valued one-form on $P \times \R^{2n}$.
The difference $\omega_0 - d \left< \Phi_{\R^{2n}}, \Theta \right> $
descends to a closed two-form $\omega_\Theta$ on $P \times_G \R^{2n}$ 
with the desired properties.

If a compact Lie group $K$ acts on $E$ by bundle automorphisms,
then, by averaging, we can choose the connection one-form $\Theta$ 
to be $K$-invariant.  The resulting two-form $\omega_\Theta$
is then $K$-invariant.
\end{proof}


Let $(B,\omega_B)$ be a symplectic manifold, let $\pi \colon E \to B$
be a symplectic vector bundle, and let a compact Lie group $K$
act on $E$ by bundle automorphisms that descend to symplectomorphisms
of $(B,\omega_B)$.   Let $\omega_\Theta$ be a closed two-form
as in Lemma~\ref{omegaTheta}.
The corresponding \emph{coupling form} is
\begin{equation} \labell{coupling}
\omega_E := \pi^* \omega_B + \omega_\Theta .
\end{equation}
It has the following properties.
\begin{enumerate}
\item
Its pullback to the fibres coincides with the fibrewise symplectic forms.
\item
Its pullback by the zero section coincides with $\omega_B$.
\item
At the points of the zero section, the fibres of $E$
are $\omega_E$-orthogonal to the zero section.
\end{enumerate}
Consequently, $\omega_E$ is non-degenerate near the zero section,
and the normal bundle of the zero section in $(E,\omega_E)$
is isomorphic to $E \to B$ as symplectic vector bundles.

\begin{Remark}
The construction that we just described is Shlomo Sternberg's
\emph{minimal coupling}~\cite{sternberg}.
For more details, see chapter 1 of the book~\cite{GLS}
by Guillemin, Lerman, and Sternberg. 
\end{Remark}

\begin{Lemma} \labell{local flow is global flow}
Let $\pi \colon E \to B$ be a symplectic vector bundle, and let 
the circle group act on $E$ by fibrewise linear symplectic transformations.  
Let $\omega_E$ be an invariant closed two-form on the total space of $E$ 
whose pullback to the fibres coincides with the fibrewise symplectic form.
Let $H \colon E \to \R$ be the fibrewise quadratic momentum map
for the fibrewise circle action.  Then $H$ is also a momentum map
for the circle action on the total space $(E,\omega_E)$.
\end{Lemma}

\begin{proof}
Let $\xi$ denote the vector field that generates the circle action.
The one-form $\iota(\xi) \omega_E$ is closed and vanishes on the zero
section, so it is exact; so there exists a momentum map
$$ H_E \colon E \to \R $$
for the circle action on $(E,\omega_E)$. 
\ ($\omega_E$ need not be nondegenerate.)
Because $\xi = 0$ along the zero section, $H_E$ is locally constant
along the zero section.  Adding a locally constant function,
we may assume that $H_E$ vanishes along the zero section.
Because the pullback of $\omega_E$ to the fibres coincides with the
fibrewise symplectic form, the restriction of $H_E$ to each fibre
is also a momentum map for the fibrewise circle action;
because it vanishes on the zero section, it must coincide
with the fibrewise quadratic momentum map,~$H$.
\end{proof}

Now let $B$ be a connected component of the fixed point set
of a circle action on a symplectic manifold $(M,\omega)$,
and let $E$ be the normal bundle of $B$ in $M$.
Then $B$ is a symplectic submanifold, and $E$ can be identified
with the symplectic orthocomplement of $TB$ in $TM|_B$.
Moreover, the circle action on $M$ induces a circle action on $E$
by bundle automorphisms.
Let $\omega_E$ be an invariant coupling form on $E$, as described above.
In this setting,
we have the following symplectic tubular neighbourhood theorem.

\begin{equation} \labell{lnf}
\begin{minipage}{0.8\textwidth} 
There exists an invariant neighbourhood $U$ of the zero section in $E$
and an equivariant symplectic open embedding
$ (U,\omega_E) \to (M,\omega) $
whose restriction to the zero section is the identity map on $B$
and, under the natural identification of $TE|_B$ with $TM|_B$,
whose differential is the identity map at every point of $B$.
\end{minipage}
\end{equation}

(This follows from the classical tubular neighbourhood theorem 
in differential topology,
combined with Theorem 4.1 of Alan Weinstein's paper \cite{weinstein0},
keeping track of a group action 
as explained in the last paragraph of Section 3 of \cite{weinstein0}.)

\bigskip

\begin{Lemma} \labell{decompose}
Let $E \to B$ be a symplectic vector bundle
with a fibrewise circle action.
Suppose that all the weights are negative;
denote the weights by $-k_1,\ldots,-k_s$
and let 
$$ k = \min \{ k_1, \ldots, k_s \} - 1 .$$
Then there exist $k+1$ commuting fibrewise circle actions on $E$,
with fibrewise quadratic momentum maps
$$ H_{(j)} \colon E \to \R \quad , \quad 0 \leq j \leq k , $$
with the following properties.
\begin{enumerate}
\item[(i)]
The $H_{j}$ are are negative outside the zero section, 
and their sum
$$ H_{(0)} + H_{(1)} + \ldots + H_{(k)} $$
is the fibrewise quadratic momentum map for the given circle action on $E$.
\item[(ii)]
There exist arbitrarily small neighbourhoods $U$ of the zero section
that are invariant under the $k+1$ commuting fibrewise circle actions,
whose intersections with the fibres of $E$ are connected,
and that satisfy
\begin{equation} \labell{inf}
\inf\limits_U \sum_{j=0}^k H_{(j)} = \sum_{j=0}^k \inf\limits_U H_{(j)}.
\end{equation}
\end{enumerate}
\end{Lemma}


\begin{Remark} \labell{rk:fibrewise}
This lemma is about fibrewise structures.  It does not require
a two-form on the total space of $E$.
\end{Remark}

\begin{proof}
Fix an invariant fibrewise compatible complex structure and Hermitian metric 
on $E$.
(Because the weights are all negative, such a metric is unique.)
After eliminating repetitions, we may assume that the 
positive numbers $k_1,\ldots,k_s$ are distinct.
Let $E_j$ denote the $-k_j$th weight space for the 
circle action on $E$, so that $E$ decomposes 
as $E_1 \oplus \ldots \oplus E_s$.
The function 
$$ z \mapsto \pi \| z_j \|^2,$$
for $z = (z_1,\ldots,z_s) \in E_1 \oplus \ldots \oplus E_s$,
is a momentum map for the scalar multiplication circle action
on the $E_j$ factor.
The quadratic fibrewise momentum map for the given circle action 
on $E$ is
\begin{equation} \labell{mm on E}
- k_1 \pi \| z_1 \|^2 - \ldots - k_s \pi \| z_s \|^2. 
\end{equation}

For each $j$, because $k_j \geq k+1$, we can decompose $k_j$ as
$$ k_j = a_{0j} + a_{1j} + \ldots + a_{kj} $$
where $a_{ij}$ are positive integers. 
We get a decomposition of~\eqref{mm on E} into the sum
$$ H_{(0)} + H_{(1)} + \ldots + H_{(k)} , $$
where
$$ H_{(i)} (z) = - a_{i1} \pi \|z_1\|^2 - \ldots - a_{is} \pi \|z_s\|^2 $$
generates the fibrewise circle action on $E$ that rotates the $j$th coordinate 
with speed $-a_{ij}$.  
Because the numbers $a_{ij}$ are positive, 
the functions $H_{(j)}$ are negative outside the zero section.

Finally, every neighbourhood of the zero section contains a neighbourhood
of the form $U = \bigcap\limits_{j=1}^s \{ z \ | \ \pi \| z_j \|^2 < \eps \}$,
and such a $U$ satisfies the requirements of part (ii) of the lemma.
\end{proof}

\begin{Lemma} \labell{vf}
Let $(B,\omega_B)$ be a symplectic manifold,
let $E \to B$ be a symplectic vector bundle,
let $\omega_{E}$ be a coupling form (see~\eqref{coupling}),
and let $U$ be an open subset of $E$
where $\omega_{E}$ is non-degenerate.
Let $k$ be a positive integer.
Suppose that there exist $k$ sections of $E$ that are nowhere all vanishing.
Then there exist smooth functions $g_1,\ldots,g_k \colon U \to \R$
whose Hamiltonian vector fields are nowhere all tangent
to the zero section of $E$.
\end{Lemma}

\begin{proof}
Let $\xi_1,\ldots,\xi_k$ be sections of $E \to B$ 
that are nowhere all vanishing.
Fix $j \in \{ 1 , \ldots , k \}$.
Let $g_j \colon E \to \R$ be the function
whose restriction to the fibre over $b$
is the linear functional $\iota(\xi_j(b)) \omega_{E_b}$,
where $\omega_{E_b}$ is the symplectic form on the fibre over $b$.
Identifying $E$ with the vertical subbundle of $TE|_B$, 
the Hamiltonian vector field of $g_j|_{U}$
is equal to $\xi_j$ at the points of the zero section.
Because $\xi_1,\ldots,\xi_k$ are vertical and nowhere all vanishing,
they are nowhere all tangent to the zero section.
\end{proof}

\begin{Remark} \labell{exist sections}
Each of the following assumptions on an oriented vector bundle $E \to B$
implies that there exist $k$ sections that are nowhere all vanishing.
Let $e(E)$ denote the Euler class of the bundle.
\begin{enumerate}
\item[(a)] $\dim B < k \rank E$.
\item[(b)] $\dim B = k \rank E$ and $e(E)^k = 0$.
\item[(c)] $\rank E = 2$, $e(E) = 0$, and $k \geq 1$.
\end{enumerate}
\end{Remark}

\begin{proof}
If $\dim B < k \rank E$ then, for any generic choice of $k$ sections of $E$,
these sections are nowhere all vanishing.  This gives (a).
If the rank of a bundle is equal to the dimension of its base,
or if the rank of a bundle is equal to $2$,
then the Euler class being zero implies the existence
of a non-vanishing section.  Applying this to the bundle
$\underbrace{E \oplus \ldots \oplus E}_{\text{$k$ times}}$
gives (b); applying it to $E$ gives~(c).
\end{proof}

\bigskip

In preparation for the proof of Theorem~\ref{shorten w subbundle},
which involves a \emph{subbundle} of the normal bundle,
we now give a refinement of the minimal coupling construction.
Let $(B,\omega_B)$ be a symplectic manifold,
let $\pi \colon E \to B$ be a symplectic vector bundle,
and let a compact Lie group $K$ act on $E$ by symplectic bundle automorphisms.
Let $E'$ be a $K$-invariant symplectic sub-bundle of $E$.

Let $E''$ be the symplectic orthocomplement of $E'$ in $E$.
Then
$E = E' \oplus E''$,
and we have a pullback diagram:
\begin{equation} \labell{pullback diagram}
\begin{CD}
E @> p'' >> E'' \\
@V p' VV @V \pi'' VV \\
E' @> \pi' >> B.
\end{CD} 
\end{equation}
Let $\omega_{\Theta'}$ and $\omega_{\Theta''}$ 
denote $K$-invariant two-forms on $E'$ and $E''$ with the properties
listed in Lemma~\ref{omegaTheta}.
The corresponding coupling forms on $E'$ and on $E$ are
$$ \omega_{E'} = {\pi'}^* \omega_B + \omega_{\Theta'}$$
and
\begin{eqnarray}
\omega_E
& = & \pi^* \omega_B + {p'}^* \omega_{\Theta'} + {p''}^* \omega_{\Theta''} 
\nonumber \\
& = & {p'}^* \omega_{E'} + {p''}^* \omega_{\Theta''} .
\labell{wha}
\end{eqnarray}
%

The proof of Theorem~\ref{shorten w subbundle}
will use the following lemma.
In it, we view $E'$ as the subset of $E$
consisting of the set of points whose fibrewise $E''$ coordinate is zero,
i.e., the preimage under $p''$ of the zero section of $E''$.

\begin{Lemma} \labell{Hamiltonians separate}
Let $(B,\omega_B)$ be a symplectic manifold,
let $E \to B$ be a symplectic vector bundle,
and let $E'$ and $E''$ be symplectic sub-bundles
such that $E = E' \oplus E''$.
Let $\omega_{E'}$ be a coupling form on $E'$ and let $\omega_E$
%
%
be a coupling form on $E$ that is compatible with $\omega_{E'}$
in the sense of~\eqref{wha}.
Let $U'$ be a neighbourhood of the zero section in $E'$
where $\omega_{E'}$ is non-degenerate, and let $U$ be a neighbourhood
of the zero section in $E$ where $\omega_E$ is non-degenerate
and that is contained  in ${p'}^{-1}(U')$.  

Let $g'_1, \ldots, g'_k \colon U' \to \R$
be smooth functions whose Hamiltonian vector fields in $U'$
are nowhere all tangent to the zero section of $E'$.
Let $g_1,\ldots,g_k \colon U \to \R$ be their pullbacks
under $p'|_U \colon U \to U'$.
Then, at each point of the zero section of $E$,
the Hamiltonian vector fields of $g_1,\ldots,g_k$
are tangent to $E'$ and are not all tangent
to the zero section of $E$.
\end{Lemma}

\begin{proof}
Viewing $E'$ as a subset of $E$
(namely, the $p''$-preimage of the zero section of $E''$),
for every point $q \in E'$ 
we have a natural decomposition
\begin{equation} \labell{natural decomp}
 TE|_q = TE'|_q \oplus E''|_{\pi(q)}.
\end{equation}
Let $\xi'_1,\ldots,\xi'_k$ be the Hamiltonian vector fields
of $g'_1,\ldots,g'_k$ in $U'$.
Let $\tilde{\xi'}_1,\ldots,\tilde{\xi'}_k$ be the corresponding
sections of $TE|_{E' \cap U}$.
That is, under the decomposition~\eqref{natural decomp},
we have 
$$ \tilde{\xi'}_j = \xi'_j \oplus 0.$$
In particular, $\tilde{\xi'}_j$ are tangent to the subset $E'$ of $E$,
and they satisfy $p'_* \tilde{\xi'}_j = \xi'_j$.
The vectors ${p''}_* \tilde{\xi'}_j$
are tangent to the zero section of $E''$,
so they are in the null space of $\omega_{\theta''}$.
It is enough to show that, along $E' \cap U$,
the Hamiltonian vector fields of $g_1,\ldots,g_k$
are equal to $\tilde{\xi'}_1, \ldots, \tilde{\xi'}_k$.
Fix $j \in \{ 1 , \ldots , k \}$.
We now compute: along $E' \cap U$,
\begin{eqnarray} 
 dg_j & = & {p'}^* dg'_j  \qquad \text{ by the definition of $g_j$} 
\nonumber \\
      & = & - {p'}^* \iota(\xi'_j) \omega_{E'}  
 \qquad \text{by the definition of $\xi'_j$ } 
\nonumber \\
      & = & - \iota(\tilde{\xi'}_j) {p'}^* \omega_{E'}  
 \qquad \text{because $p'_* \tilde{\xi'}_j = \xi'_j$ }
\nonumber \\
& = & - \iota(\tilde{\xi'}_j) 
        \left( \omega_E - {p''}^* \omega_{\Theta''}  \right)
 \qquad \text{by~\eqref{wha}} \nonumber \\
& = & - \iota(\tilde{\xi'}_j) \omega_E
 \qquad \text{because ${p''}_* \tilde{\xi'}_j$ 
          is in the null space of $\omega_{\theta''}$}. \nonumber
\end{eqnarray}
So $ dg_j = - \iota(\tilde{\xi'}_j) \omega_E $, as required.
\end{proof}


We will also use the following variant of Lemma~\ref{decompose}.

\begin{Lemma} \labell{decompose2}
Let $E \to B$ be a symplectic vector bundle
with a fibrewise circle action whose weights are all negative.
Let 
$$ E = E' \oplus E'' $$
be a decomposition of $E$ into two invariant symplectic subbundles
that are fibrewise symplectically orthogonal. 
Let $-k_1, \ldots, -k_s$ be the distinct weights for the circle action on $E'$.
Let
$$ k = \min \{ k_1, \ldots, k_s \} -1 .$$
Then there exist $k+1$ commuting fibrewise circle actions on $E$,
with fibrewise quadratic momentum maps 
$$ H_{(j)} \colon E \to \R \quad , \quad 0 \leq j \leq k ,$$
that preserve the decomposition $E = E' \oplus E''$,
that commute with the given circle action on each of $E'$ and $E''$,
and that have the following properties.

\begin{enumerate}

\item[(i)]
The function $H_{(0)}$ vanishes on the zero section on $E$
and is negative outside it,
the functions $H_{(1)}$, $\ldots$, $H_{(k)}$ all vanish on $E''$
and are negative outside it, and the sum
$$ H_{(0)} + H_{(1)} + \ldots + H_{(k)} $$
is the fibrewise quadratic momentum map for the given circle action on $E$. 

\item[(ii)]
There exist arbitrarily small neighbourhoods $U$ of the zero section in $E$
that are invariant under these $k+1$ circle actions
and under the given circle actions on $E'$ and on $E''$,
whose intersections with the fibres of $E$ are connected,
and that satisfy
\begin{equation} \labell{infs}
   \inf\limits_{U} \sum_{j=0}^k H_{(j)} 
   = \sum_{j=0}^k \inf\limits_{U} H_{(j)} .
\end{equation}

\end{enumerate}
\end{Lemma}

\begin{proof}
Fix an invariant fibrewise compatible complex structure and Hermitian metric
on each of $E'$ and $E''$.
Let $E_j$ denote the $-k_j$th weight space for the circle action on $E'$,
so that $E$ decomposes as $E_1 \oplus \ldots \oplus E_s \oplus E''$. 
The function
$$ z \mapsto \pi \| z_j \|^2 ,$$
for $z = (z_1,\ldots,z_s,z'') \in E_1 \oplus \ldots \oplus E_s \oplus E''$,
is a momentum map for the scalar multiplication circle action on the $E_j$
factor.
The quadratic fibrewise momentum map for the given circle action on $E$ is
$$ - k_1 \pi \| z_1 \|^2 - \ldots - k_s \pi \| z_s \|^2 - \hatN(z'') $$ 
where $\hatN \colon E'' \to \R$ is the negative 
of the fibrewise quadratic momentum map
for the fibrewise circle action on $E''$.
Decompose each $k_j$ as
$$ k_j = a_{0j} + a_{1j} + \ldots + a_{kj} $$
where $a_{ij}$ are positive integers.
The lemma then holds with
$$ H_{(i)}(z) = - a_{i1} \pi \| z_1 \|^2 - \ldots - a_{is} \pi \| z_s \|^2 $$
for $1 \leq i \leq k$,
with
$$ H_{(0)}(z) = - a_{01} \pi \| z_1 \|^2 - \ldots - a_{0s} \pi \| z_s \|^2 
                - \hatN(z''), $$
and with neighbourhoods of the zero section of the form
$$ U = \{ z \ | \ \pi \| z_j \|^2 < \eps \text{ for $j=1,\ldots,s$, and } 
\hatN(z'') < \eps \} $$
for small $\eps > 0$.

\end{proof}

\section{Shortening on a manifold with an arbitrary maximum}
\labell{sec:nonisolated}

In this section we prove Theorems~\ref{thm:nonisolated}
and~\ref{shorten w subbundle},
with Remarks~\ref{rk:nonisolated} and~\ref{rk:shorten w subbundle}.

We begin by recalling the statement of Theorem~\ref{thm:nonisolated}
and Remark~\ref{rk:nonisolated}:

\begin{quotation}
Let $(M,\omega)$ be a compact connected symplectic manifold
with a Hamiltonian circle action.
Let $\Bmax$ be the set where the momentum map is maximal.
Let $-k_1,\ldots,-k_s$ denote the distinct weights
for the circle action on the normal bundle of $\Bmax$,
and let $k := \min \{ k_1,\ldots,k_s \} - 1$.
Suppose that $\Bmax$ is symplectically $k$-displaceable
in every neighbourhood.
Then the circle action can be deformed through loops in $\Ham(M,\omega)$
into a loop of smaller Hofer length.

We can choose the deformation such that the positive Hofer length 
never exceeds the initial one and the negative Hofer length remains constant.
Specifically, let $H$ denote the momentum map. 
For every positive number $\alpha$,
we can choose the deformation such that the normalized deformed Hamiltonians
remain between than $H(\Bmax)-\alpha$ and $H(\Bmax)$ on the set
$\{ m \in M \ | \ H(m) > H(\Bmax) - \alpha \}$
and coincide with~$H$ outside this set.
\end{quotation}

\begin{proof}[Proof of Theorem~\ref{thm:nonisolated}
and Remark~\ref{rk:nonisolated}]

Let $E$ denote the normal bundle of $\Bmax$ in $(M,\omega)$;
recall that its fibrewise symplectic structure and circle action
are induced from $M$.
Apply Lemma~\ref{decompose} to obtain $k+1$ commuting
fibrewise circle actions on~$E$, with fibrewise quadratic momentum maps
$$ H'_{(j)} \colon E \to \R \quad , \quad 0 \leq j \leq k  $$
that are negative outside the zero section, whose sum
$$ H'_{(0)} + H'_{(1)} + \ldots + H'_{(k)} $$
is the fibrewise quadratic momentum map for the given circle action on $E$
(the one induced from $M$),
and such that part~(ii) of Lemma~\ref{decompose} holds
(with ``$H_{(j)}$" replaced by ``$H'_{(j)}$").

These circle actions fit together into a (not necessarily faithful)
action of the torus $T := (S^1)^{k+1}$,
whose diagonal acts by the given circle action on~$E$.
Apply Lemma~\ref{omegaTheta} and Equation~\eqref{coupling}, with $K=T$,
to obtain a $T$-invariant coupling form $\omega_E$ on $E$. 
Consider the given circle action on $M$ and the induced circle action on $E$,
and apply the symplectic tubular neighbourhood theorem~\eqref{lnf},
to obtain an invariant neighbourhood $U$ of the zero section in $E$
on which $\omega_E$ is nondegenerate
and an equivariant open symplectic embedding
$$ i \colon (U,\omega_E) \to (M,\omega) $$ 
whose restriction to the zero section is the identity map on $\Bmax$
and, under the natural identification of $TE|_{\Bmax}$ with $TM|_{\Bmax}$,
whose differential is the identity map at every point of~$\Bmax$.

Fix a positive number, $\alpha > 0$.
After possibly shrinking $U$, we may assume that
$i(U) \subset \{ m \in M \ | \ H(m) > H(\Bmax) - \alpha \}$,
and that $U$ satisfies the requirements in part~(ii) of Lemma~\ref{decompose}:
\ $U$ is invariant under the entire $T$ action on $E$, 
the intersections of $U$ with the fibres of $E$ are connected, and 
$$ \inf\limits_U \sum_{j=0}^k H'_{(j)}
 = \sum_{j=0}^k \inf\limits_U H'_{(j)} .$$

Let $ H \colon M \to \R $ be the momentum map on $M$.
We have 
\begin{align} \labell{left right H}
 i^*H & = \left. \left( 
      H(\Bmax) + H'_{(0)} + H'_{(1)} + \ldots + H'_{(k)} \right) \right|_{U} \\
 & = H_{(0)} + H_{(1)} + \ldots + H_{(k)}, \nonumber
\end{align}
where $H_{(0)} = H(\Bmax) + H'_{(0)}|_{U}$
and where $H_{(j)} = H'_{(j)}|_{U}$ for each $1 \leq j \leq k$.
(The intersection with $U$ of each fibre of $E$ is an invariant
connected open neighbourhood of the origin in the fibre.  
In this neighbourhood,
the left and right hand sides of~\eqref{left right H}
are momentum maps for the same circle action with the same value at the origin,
so they are equal.)

Because the $T$ action on $E$ preserves $\omega_E$ 
and by Lemma~\ref{local flow is global flow},
the maps $H'_{(0)}$, $H'_{(1)}$, $\ldots$, $H'_{(k)}$
are momentum maps for the $k+1$ commuting circle actions
on the total space $(E,\omega_E)$ and not only fibrewise.
So $H_{(0)}, H_{(1)}, \ldots , H_{(k)}$ 
generate commuting circle actions on $(U,\omega_E)$.

We are now in the setup of Lemma~\ref{lem polterov k}.
The sets $\Bmax^{(j)}$ where the functions $H_{(j)}$ attain
their maximum are all equal to the zero section of $E$.
This implies the condition~\eqref{included} 
of Lemma~\ref{lem polterov k}.
It also reduces the last condition of Lemma~\ref{lem polterov k}
to the condition that the zero section is symplectically $k$-displaceable 
in $U$.  This condition then follows from that $\Bmax$ is symplectically 
$k$-displaceable in $i(U)$.  Theorem~\ref{thm:nonisolated} 
and Remark~\ref{rk:nonisolated}
then follows from Lemma~\ref{lem polterov k}.

\end{proof}

Next, we prove Theorem~\ref{shorten w subbundle}
and Remark~\ref{rk:shorten w subbundle}.
We recall the statement:

\begin{quotation}
Let $(M,\omega)$ be a compact connected symplectic manifold
with a Hamiltonian circle action.
Let $\Bmax$ be the set where the momentum map is maximal.
Let $E' \to \Bmax$ be an $S^1$-invariant subbundle
of the normal bundle to $\Bmax$ in $M$,
and let $-k'_1, \ldots, -k'_{s'}$ be the distinct weights
for the circle action on $E'$.
Let 
$$ k' = \min \{ k'_1, \ldots, k'_{s'} \} -1 . $$
Suppose that 
\begin{equation}  
 \text{$E'$ has $k'$ sections that are nowhere all vanishing.}
\end{equation}
Then the circle action can be deformed through loops in $\Ham(M,\omega)$
into a loop of smaller Hofer length.

We can choose the deformation such that the positive Hofer length 
never exceeds the initial one and the negative Hofer length remains constant.
Specifically, let $H$ denote the momentum map. 
For every positive number $\alpha$,
we can choose the deformation such that the normalized deformed Hamiltonians
remain between $H(\Bmax)-\alpha$ and $H(\Bmax)$ on the set
$\{ m \in M \ | \ H(m) > H(\Bmax) - \alpha \}$
and coincide with $H$ outside this set.
\end{quotation}

{

\begin{proof}[Proof of Theorem \ref{shorten w subbundle}
and Remark~\ref{rk:shorten w subbundle}]

Let $E$ denote the normal bundle of $\Bmax$ in $(M,\omega)$;
recall that its fibrewise symplectic structure and circle action
are induced from $M$.
Let $E''$ be the fibrewise symplectic orthocomplement of $E'$ in $E$,
so $E = E' \oplus E''$.
Apply Lemma~\ref{decompose2} to obtain $k'+1$ commuting fibrewise 
circle actions on $E$,
with fibrewise quadratic momentum maps
$$ H'_{(j)} \colon E \to \R \quad , \quad 0 \leq j \leq k' $$
such that $H'_{(0)}$ vanishes on the zero section of $E$
and is negative outside it,
$H'_{(1)}, \ldots, H'_{(k')}$ vanish on $E''$ and are negative outside it,
the sum $H'_{(0)} + H'_{(1)} + \ldots + H'_{(k')}$
is the fibrewise quadratic momentum map for the given circle action on $E$,
and part (ii) of Lemma~\ref{decompose2} holds
(with ``$k$" and ``$H_{(j)}$" replaced by ``$k'$" and ``$H'_{(j)}$").

These $k'+1$ circle actions, together with the given circle actions
on $E'$ and on $E''$, 
%
%
fit together into a fibrewise action of a torus $T$. 
Apply Lemma~\ref{omegaTheta} and Equations~\eqref{coupling}
and~\eqref{wha} to obtain $T$ invariant coupling forms 
$\omega_{E'}$ on $E'$ and $\omega_E$ on $E$.
Consider the given circle action on $M$ and the induced circle action on~$E$,
and apply the symplectic tubular neighbourhood theorem~\eqref{lnf}
to obtain an invariant neighbourhood $U$ of the zero section in $E$
on which $\omega_E$ is nondegenerate
and an equivariant symplectic open embedding
$$ i \colon (U,\omega_E) \to (M,\omega) $$
whose restriction to the zero section is the identity map on $\Bmax$
and, under the natural identification of $TE|_{\Bmax}$ with $TM|_{\Bmax}$, 
whose differential is the identity map at every point of~$\Bmax$.

Fix a positive number, $\alpha > 0$.
After possibly shrinking $U$, we may assume that 
$i(U) \subset \{m \in M \ | \ H(m) > H(\Bmax)-\alpha \}$
and that $U$ satisfies the requirements in part (ii) 
of Lemma~\ref{decompose2}:
\ $U$ is $T$ invariant,
the intersections of $U$ with the fibres of $E$ are connected,
and
$$ \inf\limits_U \sum_{j=0}^{k'} H'_{(j)}
 = \sum_{j=0}^{k'} \inf\limits_U H'_{(j)} . $$

Let $H \colon M \to \R$ denote the momentum map on $M$.  We have
\begin{align} \labell{left right H again}
 i^*H & = \left. \left( H(\Bmax) + H'_{(0)} + \ldots + H'_{(k')}  
                 \right)\right|_{U} \\
\nonumber &  = H_{(0)} + H_{(1)} + \ldots + H_{(k')},
\end{align}
where 
$$ H_{(0)} = \left.\left( H(\Bmax) + H'_{(0)} \right)\right|_U $$
and where 
$$ H_{(j)} = \left.\left( H'_{(j)} \right)\right|_U 
             \quad \text{ for each } \quad 1 \leq j \leq k'.$$
(The intersection of $U$ with each fibre of $E$
is an invariant connected open neighbourhood of the origin in the fibre.
In this neighbourhood, the left and right hand sides 
of~\eqref{left right H again}
are momentum maps for the same circle action
with the same value at the origin, so they are equal.)

Because the $T$ action on $E$ preserves $\omega_E$ 
and by Lemma~\ref{local flow is global flow}, the maps
$H'_{(0)}$, $H'_{(1)}$, $\ldots$, $H'_{(k')}$
are momentum maps for the $k'+1$ commuting circle actions
on the total space $(E,\omega_E)$ and not only fibrewise.
So $H_{(0)}$, $H_{(1)}$, $\ldots$, $H_{(k')}$
generate commuting circle actions on $(U,\omega_E)$.

We are now in the setup of Lemma~\ref{lem polterov k}.
The set $\Bmax^{(0)}$ where $H_{(0)}$ attains its maximum
is equal to the zero section, and, for $1 \leq j \leq k'$,
the set $\Bmax^{(j)}$ where $H_{(j)}$ attains its maximum
is the intersection of $U$ with the image of $E''$ in $E$.
Slightly abusing notation,
we write these sets as $\Bmax$ and as $E'' \cap U$.
To apply Lemma~\ref{lem polterov k}, we need to find
symplectomorphisms $b_1,\ldots,b_{k'}$ of~$U$, each connected to the identity
through a path of compactly supported symplectomorphisms of~$U$, 
such that
\begin{equation}\labell{need}
 \Bmax \cap b_1(E'' \cap U) \cap \ldots \cap b_{k'}(E'' \cap U) = \emptyset.
\end{equation}

Let $U'$ be a neighbourhood of the zero section in $E'$
where $\omega_{E'}$ is nondegenerate.
By Lemma~\ref{vf}, there exist functions $g'_1,\ldots,g'_{k'}$ on $U'$
whose Hamiltonian vector fields are nowhere all tangent to the
zero section of~$E'$.
For $j=1,\ldots,k'$, let $g_j = \rho (g'_j \circ p')$
where $p' \colon E \to E'$ is the projection map 
and where $\rho \colon U \to \R$ 
is a cutoff function that is equal to~$1$ near the zero section of $E$
and has compact support in $U$.
Lemma~\ref{Hamiltonians separate} implies that,
at each point of $\Bmax$,
the Hamiltonian vector fields of $g_1,\ldots,g_{k'}$
are tangent to~$E'$ and are not all tangent to $\Bmax$.
Let $\psi^{g_1}_t, \ldots, \psi^{g_{k'}}_t$ be the corresponding flows.
If $t$ is a sufficiently small positive number, the intersection
$\Bmax \cap \psi^{g_1}_t(E'' \cap U) \cap \ldots 
       \cap \psi^{g_{k'}}_t(E'' \cap U)$
is empty.
Lemma~\ref{lem polterov k}, 
with $b_1,\ldots,b_{k'}$ taken to be $\psi^{g_1}_t, \ldots, \psi^{g_{k'}}_t$,
gives the result of Theorem~\ref{shorten w subbundle}
and Remark~\ref{rk:shorten w subbundle}.
\end{proof}

}

\section{The index of the positive Hofer length functional}
\labell{sec:savelyev}

In this section we prove Theorem \ref{thm:savelyev} 
and Remark~\ref{rk:savelyev}.  We recall the statement,
with slight change of notation ($H_M$ instead of $H$):

\begin{quotation}
Let $(M,\omega)$ be a $2n$ dimensional
compact connected manifold with a Hamiltonian
circle action.  
Suppose that the momentum map attains its maximum at an isolated fixed point.
Let $-k_1,\ldots,-k_n$ be the isotropy weights at the maximum,
with possible repetitions.
Then there exists a neighbourhood $D$ of the origin 
in $\R^{\sum 2(k_j - 1)}$, and, for each $\lambda \in D$, a loop
$\{ \psi^{(\lambda)}_t \}_{0 \leq t \leq 1}$ in $\Ham(M,\omega)$,
such that the following properties hold.
For $\lambda=0$, the loop $\{\psi^{(0)}_t\}_{0 \leq t \leq 1}$ 
is the given circle action.
The function $\lambda \mapsto \length( \{ \psi^{(\lambda)}_t \} )$ is smooth, 
$\lambda=0$ is a critical point of this function,
and the Hessian of this function at $\lambda=0$ is negative definite.

In fact, we can choose the deformation such that the negative 
Hofer length remains constant.  Specifically, let $H_M$ denote
the momentum map, and let $\pmax$ denote the point where $H_M$ attains
its maximum. For every positive number $\alpha$,
we can choose the deformation such that the normalized deformed Hamiltonian
remains between $H_M(\pmax)-\alpha$ and $H_M(\pmax)$ on the set
$\{ m \in M \ | \ H_M(m) > H_M(\pmax)-\alpha \}$
and coincides with $H_M$ outside this set.
\end{quotation}

\begin{proof}[Proof of Theorem~\ref{thm:savelyev}
and Remark~\ref{rk:savelyev}]

We set the following notation.  
$$ \phi_t \colon \C \to \C \quad , \quad \phi_t(z) = e^{-2\pi i t} z 
   \quad \text{ for } t \in \R.$$
$$ \beta_\lambda \colon \C \to \C \quad , \quad 
    \beta_\lambda(z) = z + \lambda 
   \quad \text{ for } \lambda \in \C.$$

Let $\pmax$ denote the point where the momentum map $H_M$ attains its maximum.
By the equivariant Darboux theorem, near $\pmax$, 
we can identify the manifold with $\C^n$, the action with the product action
\begin{equation} \labell{product action}
(\phi_{k_1 t}) \times \ldots \times (\phi_{k_n t}),
\end{equation}
and the momentum map with the function 
$$ z \mapsto H_M(\pmax) - N(z) ,$$
where 
\begin{equation} \labell{original N}
 N(z) = \pi \left( k_1 |z_1|^2 + \ldots + k_n |z_n|^2 \right) .
\end{equation}

(We will use variants of the symbol $N$
to denote variations and combinations of norm-square functions,
and we will use variants of the symbol $H$ to denote momentum maps.)

If $k$ is an integer greater than $1$,
the non-effective circle action $\phi_{kt}$
on $\C$ has the deformation
\begin{equation} \labell{deformation on C}
\phibar_t^{(\lambda)} \, := \,
   \phi_t \circ \beta_{\lambda_{k-1}} \circ \phi_t
   \circ \ldots 
   \circ \beta_{\lambda_{2}} \circ \phi_t 
   \circ \beta_{\lambda_{1}} \circ \phi_t 
   \circ \beta_{-\lambda_1-\ldots-\lambda_{k-1}} ,
\end{equation}
parametrized by $(\lambda) = (\lambda_1,\ldots,\lambda_{k-1}) \in \C^{k-1}$
and generated by 
\begin{multline} \labell{sum of terms}
\Hbar^{(\lambda)}_t :=  
 H + H \beta_{\lambda_{k-1}}\inv \phi_t\inv
   + H \beta_{\lambda_{k-2}}\inv \phi_t\inv \beta_{\lambda_{k-1}}\inv \phi_t\inv
   + \ldots \\
 \ldots
   + H \beta_{\lambda_{1}}\inv \phi_t\inv \beta_{\lambda_{2}}\inv \phi_t\inv
       \cdots \beta_{\lambda_{k-1}}\inv \phi_t\inv ,
\end{multline}
where $H(z) = -\pi |z|^2$.
(We will use the symbol $(\lambda)$, with brackets, for a parameter 
in $\C^{k-1}$ or in $\C^{\sum (k_i-1)}$
and the symbol $\lambda$, without brackets, for a parameter in $\C$.)

Applying such a deformation to each factor of \eqref{product action},
we get a family of deformations of the product action.
This family is parametrized by elements 
of $\C^{\sum (k_i-1)}$.
When we view the parameter $(\lambda)$ 
as an element $((\lambda)_{1},\ldots,(\lambda)_{n})$
of the $n$-fold product $\C^{k_1-1} \times \ldots \times \C^{k_n-1}$,
the deformation is
\begin{equation} \labell{product}
\phibar^{(\lambda)_{1}} \times
 \ldots \times \phibar^{(\lambda)_{n}},
\end{equation}
and it is generated by 
\begin{equation} \labell{deformation}
 H^{(\lambda)}_t(z) := \Hbar^{(\lambda)_{1}}_t(z_1) + 
 \ldots + \Hbar^{(\lambda)_{n}}_t(z_n).
\end{equation}

With the help of cut-off functions,
this deformation can be plugged into the original manifold
to obtain a family of deformations of the original circle action,
parametrized by $\lambda$s
in a neighbourhood $D$ of the origin in~$\C^{\sum (k_i-1)}$,
and generated by a function that in a neighbourhood of $p$
can be identified with $z \mapsto H_M(\pmax) + H^{(\lambda)}_t(z)$
and that attains its maximum in that neighbourhood.
We give the details later.

We now turn to the relevant computation on $\C^n$.
We consider the function
$$ (\lambda) \mapsto \int_0^1 \left( 
    \max\limits_{z \in \C^n} H^{(\lambda)}_t(z) \right) dt $$
from $\C^{\sum(k_i-1)}$ to $\R$.
By~\eqref{sum of terms} and~\eqref{deformation},
this function is everywhere non-positive and vanishes at the origin.  
We will show that this function is quadratic
in the real and imaginary parts of $\lambda$,
and that it is strictly negative when $\lambda \neq 0$.
This will imply that the function is smooth and that its Hessian 
at the origin is negative definite, as required.

By~\eqref{deformation}, it is enough to consider separately
each of the functions
$(\lambda)_{i} \mapsto 
 \int_0^1 \max\limits_\C \Hbar^{(\lambda)_{i}}_t dt$.
Omitting the index $i$, for each $0 \leq t \leq 1$ we consider the function
\begin{equation} \labell{need to show}
(\lambda) \mapsto \int_0^1 \left( \max_{z \in \C} 
                                  \Hbar^{(\lambda)}_t \right) dt
\end{equation}
from $\C^{k-1}$ to $\R$.  
We claim that it is a quadratic function of the real and imaginary parts
of $\lambda$ whose coefficients are trigonometric polynomials in $t$,
and that it is everywhere non-positive, vanishes at $\lambda = 0$,
and strictly negative for $\lambda \neq 0$ when $t=0$.

Expanding the expression~\eqref{sum of terms}, we get
that $\Hbar^{(\lambda)}_t(z)$ is equal to $-\pi$ times  
\begin{multline} \labell{crawford}
  |z|^2 + |e^{2 \pi i t} z - \lambda_{k-1} |^2 
      + | e^{2 \pi i \cdot 2t} z
          - e^{2 \pi i \cdot t} \lambda_{k-1} - \lambda_{k-2}  |^2 
    + \ldots \\
    + | e^{2\pi i (k-1)t}z - e^{2\pi i (k-2) t}\lambda_{k-1} 
    - \ldots
    - e^{2 \pi i \cdot 2t}\lambda_3 - e^{2 \pi i t} \lambda_2 - \lambda_1|^2.
\end{multline}
Expanding further, this is equal to 
\begin{multline} \nonumber
 |z|^2 + |z|^2 + |\lambda_{k-1}|^2 
                      - 2\real(e^{2 \pi i t} z \cdot \lambda_{k-1}) \\
 + |z|^2 + |e^{2 \pi i t} \lambda_{k-1} + \lambda_{k-2}|^2 
               - 2 \real(e^{2 \pi i \cdot 2t} z 
                   \cdot (e^{2 \pi i t} \lambda_{k-1} + \lambda_{k-2})) \\
 + \ldots \\
 + |z|^2 + | e^{2 \pi i (k-2) t} \lambda_{k-1} + \ldots +
           e^{2 \pi i \cdot 2t} \lambda_3 + e^{2 \pi i t} \lambda_2 
                                                          + \lambda_1|^2 \\
       - 2 \real ( e^{2 \pi i (k-1) t} z \cdot (
                    e^{2 \pi i (k-2) t} \lambda_{k-1}
                  + \ldots 
                  + e^{2 \pi i \cdot 2 t} \lambda_3
                  + e^{2 \pi i t} \lambda_2 + \lambda_1) ).
\end{multline}
Thus, $\Hbar^{(\lambda)}_t(x+iy)$ 
is equal to $-\pi$ times $kx^2 + ky^2 + ax + by + c$,
where each of $a,b$ is a linear function of the real and imaginary parts
of $\lambda_1,\ldots,\lambda_{k-1}$
whose coefficients are trigonometric polynomials in $t$
and where $c$ is a quadratic function of the real and imaginary parts
of $\lambda_1,\ldots,\lambda_{k-1}$ whose coefficients are trigonometric
polynomials in $t$.
Completing the square, this function attains its minimum
when $x=-\frac{a}{2k}$ and $y = -\frac{b}{2k}$, and its minimal value is
$c - \frac{a^2}{4k} - \frac{b^2}{4k}$ 
which, as required, is a quadratic function
of the real and imaginary parts of $\lambda_1,\ldots,\lambda_{k-1}$  
whose coefficients are trigonometric polynomials in $t$.

Finally, if $\lambda \neq 0$, then at least one of the summands
in~\eqref{crawford} is strictly positive when $z=0$,
whereas the first summand is strictly positive when $z \neq 0$,
so the minimal value of the function~\eqref{crawford} is strictly positive.

\bigskip

This completes our computation in $\C^n$.
We now return to our manifold $M$.
Fix a symplectomorphism $f \colon U \to M$ 
from a neighbourhood $U$ of the origin in $\C^n$
onto a neighbourhood of $\pmax$ in $M$ such that $f(0) = \pmax$
and such that $(f^* H_M)(z) = H_M(\pmax) - N(z)$,
where $N(\cdot)$ is given in~\eqref{original N}.
Given a positive number $\alpha$, we may choose $U$ such that 
$i(U) \subset \{m\in M \ | \ H_M(m) > H_M(\pmax) -\alpha \}$.

Let $\epsbar$ be a positive number that is sufficiently small
so that the set $\{ z \in \C^n \ | \ N(z) \leq \epsbar \}$
is contained in $U$.
Then the subset of $M$ where $H_M$ is $\epsbar$-close to its maximum $H_M(p)$ 
is contained in $f(U)$.
(See the facts about Hamiltonian circle actions in \S\ref{sec:preliminaries}.)

Let $\rho \colon [0,\infty) \to [0,1]$ be a smooth function
that is equal to one on $[0,\epsbar/2]$ 
and is equal to zero outside $[0,\epsbar)$.
Use the same symbol, $\rho$, to denote the function
$z \mapsto \rho(N(z))$ on $\C^n$.
Then the function $H^{(\rho\cdot \lambda)}_t(z)$ 
is equal to $H^{(\lambda)}_t(z)$ when $N(z) \leq \epsbar/2$ 
and is equal to $-N(z)$ when $N(z) > \epsbar$.

The function on $M$
\begin{equation} \labell{deformed}
{H_M}^{(\lambda)}_t := 
 \begin{cases}
H_M(\pmax) + H^{(\rho \cdot \lambda)}_t(\cdot) \circ f\inv \\
\phantom{H_M} \qquad \qquad 
 \text{ when } H_M \text{ is $\epsbar$-close to its maximum } \\
H_M \qquad \qquad \text{ elsewhere }
\end{cases}
\end{equation}
is smooth and it generates a deformation of our loop.

\smallskip

If the original momentum map $H_M$ is normalized, 
so is the deformed Hamiltonian~\eqref{deformed}.
The argument is similar to~\eqref{3}.

\smallskip

Let $D$ be a neighbourhood of the origin in $\C^{\sum(k_i-1)}$ 
that is small enough so that, for every $\lambda \in D$,
the following facts are true.
\begin{itemize}
\item[-]
For every $t$, the function $H^{(\lambda)}_t(z)$ attains its maximum
at a unique point of $\C^n$, and this point
lies in $\{ z \ | \ N(z) \leq \epsbar/4 \}$.
\item[-]
This maximum is greater than $-\epsbar/4$.
\item[-]
For every $\nu$ and $z$ such that $0 \leq \nu \leq 1$
and $\epsbar/2 \leq N(z) \leq \epsbar$,
we have $H^{(\nu \cdot \lambda)}_t(z) \leq -\epsbar/4$
for all $t$.
\end{itemize}
This is possible because $H^{(\lambda)}_t(z)$ is smooth 
in the variables $\lambda$, $z$, and $t$, is quadratic in $z$,
and is equal to $-N(z)$ when $\lambda=0$.

Then, for every $\lambda$ and $t$,
the functions $H^{(\rho \cdot \lambda)}_t(z)$ and $H^{(\lambda)}_t(z)$
have the same maximal value,
and they attain it at the same (unique) point $z$.

\smallskip

Denote by $\psi^{(\lambda)}_t$ the family of loops in $\Ham(M,\omega)$
that are generated by the functions~\eqref{deformed}.
By our choice of $D$, for every $\lambda \in D$,
the positive Hofer length of the corresponding loop is
$$ \ell_+ (\{\psi^{(\lambda)}_t\}) = 
 H_M(\pmax) + \int_0^1 \max\limits_{z \in \C^n} H^{(\lambda)}_t(z) dt.$$
By this and our computation in $\C^n$, 
the function
$$ \lambda \mapsto \ell_+( \{ \psi^{(\lambda)}_t \} ) $$
from $D$ to $\R$
takes its maximal value at the origin
and its Hessian at the origin is negative definite.
This completes the proof of Theorem~\ref{thm:savelyev} 
and Remark~\ref{rk:savelyev}.
\end{proof}


\begin{thebibliography}{HZ}
\bibitem{At}
M.~Atiyah, \emph{Convexity and commuting hamiltonians},
Bull.\ London Math.\ Soc.\ \textbf{14} (1982), 1--15.

\bibitem{Ba}
A.~Banyaga, \emph{Sur la structure du groupe des difféomorphismes 
qui préservent une forme symplectique},
Comment.\ Math.\ Helv. \textbf{53} (1978), 174--227.

\bibitem{BP}
M.~Bialy and L.~Polterovich, \emph{Geodesics of Hofer's metric
on the group of Hamiltonian diffeomorphisms},
Duke Math.\ J.\ \textbf{76} (1994), 273--292.

\bibitem{GLS}
V.\ Guillemin, E.\ Lerman, and S.\ Sternberg,
\emph{Symplectic fibrations and multiplicity diagrams},
Cambridge University Press, 1996.

\bibitem{GS:convexity}
V.~Guillemin and S.~Sternberg, \emph{Convexity properties 
of the moment mapping}, Invent.\ Math.\ \textbf{67} (1982), 491--513.

\bibitem{H:topological}
H. Hofer, \emph{On the topological properties of symplectic maps},
Proceedings of the Royal Society of Edinburgh \text{115} A (1990), 25--28.

\bibitem{H}
H.~Hofer, \emph{Estimates for the energy of a symplectic map,}
Comment.\ Math.\ Helv.\ \textbf{68} (1993), 48--72.

\bibitem{HZ}
H.~Hofer and E.~Zehnder, \emph{Symplectic Invariants and
Hamiltonian Dynamics,} Birkhauser, Boston, MA (1994).

\bibitem{KT}
Y.~Karshon and S.~Tolman,
\emph{The Gromov width of complex Grassmannians,}
Algebraic and Geometric Topology \textbf{5} (2005), 911--922.

\bibitem{LM}
F.~Lalonde and D.~McDuff, \emph{Hofer's $L^{\infty}$-geometry: 
energy and stability of Hamiltonian flows parts I and II}, 
Invent. Math. \textbf{122} (1995), 1--33 and 35--69.

\bibitem{LM:local}
F. Lalonde and D. McDuff, \emph{Local non-squeezing theorems
and stability}, Geom. and Func. Analysis vol.~\textbf{5}, 
No.~2 (1995), 364--386.

\bibitem{long}
Yi Ming Long, \emph{Geodesics in the compactly supported Hamiltnoian 
diffeomorphisms group}, Math.\ Z.\ \textbf{220} (1995), no.~2, 279--294.  

\bibitem{McD:variants}
D. McDuff, \emph{Geometric variants of the Hofer norm},
J. Sympl. Geom. \textbf{1} (2002), p.197--252.

\bibitem{MSl} 
D.~McDuff and J.~Slimowitz,
\emph{Hofer-Zehnder capacity and length minimizing Hamiltonian paths},
Geom.\ Topol.\ \textbf{5}  (2001), 799--830.

\bibitem{MT}
D. McDuff and S. Tolman,
\emph{Topological properties of Hamiltonian circle actions},
IMRP Int.\ Math.\ Res.\ Pap., (2006), 1--77.

\bibitem{MT2}
D. McDuff and S. Tolman,
\emph{On nearly semifree circle actions}, 
arXiv:math/0503467 [math.SG].

\bibitem{milnor}
J.\ Milnor, \emph{Morse theory}, Princeton University Press, 1969.


\bibitem{Pol}
L.~Polterovich, \emph{The geometry of the group of symplectic 
diffeomorphisms}, Lectures in Mathematics ETH Zürich
Birkh\"auser Verlag, Basel, 2001.

\bibitem{sav1}
Y.~Savelyev, \emph{Quantum characteristic classes
and the Hofer metric}, 
Geom.\ Topol.\ \textbf{12} (2008), no.~4, 2277--2326.

\bibitem{sav2}
Y.~Savelyev,
\emph{Virtual Morse theory on $\Omega \Ham(M,\omega)$},
J.\ Diff.\ Geom.\ \textbf{84} (2010), no.~2, 409--425.

\bibitem{sav3}
Y.~Savelyev, 
\emph{Proof of the index conjecture in Hofer geometry},
arXiv:1204.3098 [math.SG],
to appear in Math.\ Res.\ Letters.

\bibitem{Si}
K.~Siburg, \emph{New minimal geodesics in the group of symplectic 
diffeomorphisms}, Calc.\ Var.\ Part.\ Diff.\ Eq.\ \textbf{3} (1995), 
299--309.

\bibitem{sternberg}
S.\ Sternberg, \emph{Minimal coupling and the symplectic mechanics
of a classical particle in the presence of a Yang--Mills field},
Proc.\ Nat.\ Acad.\ Sci.\ U.S.A.\ \textbf{74} (1977), no.~23, 5253--5254.

\bibitem{Us}
I.~Ustilovsky, \emph{Conjugate points on geodesics of Hofer's metric}, 
Differential Geometry and its Applications \textbf{6} (1996), 327--342.

\bibitem{weinstein0}
A. Weinstein, \emph{Symplectic manifolds and their Lagrangian submanifolds},
Advances in Mathematics \textbf{6}, 329--346 (1971).

\bibitem{weinstein1}
A.\ Weinstein, \emph{Lectures on symplectic manifolds},
CBMS Lecture Notes \textbf{29}, A.M.S., Providence RI, 1977.


\end{thebibliography}
\end{document}